\documentclass[10pt,reqno]{amsart}
\usepackage{amssymb,amscd}
\newtheorem{theorem}{Theorem}[section]                    
\newtheorem{proposition}[theorem]{Proposition}            
\newtheorem{corollary}[theorem]{Corollary}                
\newtheorem{lemma}[theorem]{Lemma}
\newtheorem{remark}[theorem]{Remark} 
\newtheorem{remarks}[theorem]{Remarks}
\newtheorem{fact}[theorem]{Fact}

\newtheorem{example}[theorem]{Example}

\newtheorem{definition}{Definition}[section]
\newtheorem{assumption}[definition]{Assumption}

\begin{document}

\title[HNN Extensions and Amalgamated Free Products]{A relationship between HNN Extensions and Amalgamated Free Products in Operator Algebras$^{**}$}
\author[Y. Ueda]{Yoshimichi UEDA$^*$}
\address{Graduate School of Mathematics, Kyushu University\\ Fukuoka, 810-8560, Japan}
\email{ueda@math.kyushu-u.ac.jp}
\thanks{$^*\,$Supported in part by Grant-in-Aid for Young Scientists (B)17740096.}
\thanks{$^{**}$The revised and expanded version of this paper will be appeared in Illinois J.~Math.~under the new title ``Remarks on HNN extensions in operator algebras". This is still an old version with a few corrections. }

\maketitle

\begin{abstract} 
With a minor change made in the previous construction in \cite{Ueda:JFA05} we observe that any reduced HNN extension is precisely a compressed algebra of a certain reduced amalgamated free product in both the von Neumann algebra and the $C^*$-algebra settings. It is also pointed out that the same fact holds true even for universal HNN extensions of C$^*$-algebras. We apply the observation to the questions of factoriality and type classification of HNN extensions of von Neumann algebras and also those of simplicity and $K$-theory of (reduced or universal) HNN extensions of $C^*$-algebras.    
\end{abstract} 

\allowdisplaybreaks{

\section*{Introduction} 

This paper is a continuation of the previous work \cite{Ueda:JFA05}, where we introduced the notion of reduced HNN extensions in the category of von Neumann algebras as well as that of $C^*$-algebras, which includes, as expected, the group von Neumann algebras or the reduced group $C^*$-algebras associated with HNN extensions of groups. Firstly, with a minor change of the previous construction in \cite{Ueda:JFA05} we see that any reduced HNN extension is precisely a compressed algebra of a certain reduced amalgamated free product. The same fact is also pointed out to be valid even for universal HNN extensions of $C^*$-algebras. The observation is new even for the group von Neumann algebras and the (both reduced and universal) group $C^*$-algebras associated with HNN extensions of groups, and indeed it seems that there is no explicit counterpart in the framework of group theory. However, a similar observation was already pointed out by D.~Gaboriau \cite{Gaboriau:InventMath00} (also F.~Paulin \cite{Paulin:MarkovProcess99}) for equivalence relations. We missed it when we did \cite{Ueda:JFA05}, and comparing it with our construction of reduced HNN extensions is a starting point of the present work. Based on the observation we obtain several results on HNN extensions of von Neumann algebras and/or $C^*$-algebras: We first improve the previous factoriality and type-classification results in \cite{Ueda:JFA05}, which lead to a satisfactory answer to the questions of factoriality and type classification of HNN extensions of von Neumann algebras, say $N\bigstar_D\theta$, when both $D$ and $\theta(D)$ are (not necessarily inner conjugate) Cartan subalgebras in $N$. Here, we note that the inner conjugate case was already treated in \cite[Remark 3.7 (1)]{Ueda:JFA05} based on its particularity, and one should remind that all Cartan subalgebras in a fixed von Neumann algebra are isomorphic, which allows to take an HNN extension by a bijective $*$-homomorphism between those. We also consider the questions of simplicity and $K$-theory of (reduced or universal) HNN extensions of $C^*$-algebras. Some of the consequences here for $C^*$-algebras  should be read as improvements of previous arguments for the group $C^*$-algebras associated with HNN extensions of groups, but some others are new. 

The organization of this paper is as follows. The next \S1 is a preliminary section, where we briefly recall some of the materials in \cite{Ueda:JFA05}. In \S2 our observation mentioned above is given in the von Neumann algebra, the reduced $C^*$-algebra and the universal $C^*$-algebra settings, respectively. We also compare it with Gaboriau's observation on HNN equivalence relations, and then motivated by Gaboriau's one we show that any amalgamated free product can be described by means of a certain HNN extension by a {\it partial} $*$-isomorphism. In \S3 we give the above-mentioned applications of our observation. An appendix is presented, where we give one more factoriality result for HNN extensions of von Neumann algebras.  

\section{Reduced HNN extensions} 

Throughout this paper we follow the notational conventions in \cite{Ueda:JFA05}, which are summarized here for the reader's convenience. 

\subsection{von Neumann Algebra Setup} Let $N \supseteq D$ be $\sigma$-finite von Neumann algebras and $\theta : D \rightarrow N$ be an injective normal unital $*$-homomorphism. Assume that there are faithful normal conditional expectations $E_D^N : N \rightarrow D$, $E_{\theta(D)}^N : N \rightarrow \theta(D)$. The {\it HNN extension} of {\it base algebra} $N$ by $\theta$ with respect to $E_D^N$, $E_{\theta(D)}^N$ is a unique triple $\left(M,E_N^M : M \rightarrow N,u(\theta)\right)$ of a von Neumann algebra containing $N$, a faithful normal conditional expectation and a unitary in $M$ (called the {\it stable unitary}), which can be characterized by the following conditions: 
\begin{itemize} 
\item[(A)] 
$u(\theta)\theta(d)u(\theta)^* = d$ for every $d \in D$; 
\item[(M)] 
$E_N^M(w) = 0$ for every reduced word $w$ in $N$ and $u(\theta)$,   
\end{itemize} 
where a given word $w = u(\theta)^{\varepsilon_0} n_1 u(\theta)^{\varepsilon_1} n_2 
\cdots n_{\ell} u(\theta_{\ell})^{\varepsilon_{\ell}}$ in $N$ and $u(\theta)$ (with $n_1,\dots, n_{\ell} \in N$, $\varepsilon_0,\dots, \varepsilon_{\ell} \in \left\{1, -1\right\}$) is said to be a reduced one (or to be of reduced form) if $\varepsilon_{j-1} \neq \varepsilon_j$ implies that  
\begin{itemize}
\item[] $n_j \in N_{\theta}^{\circ} \overset{\mathrm{def}}{:=} 
\mathrm{Ker}E_{\theta(D)}^N$ when $\varepsilon_{j-1}  = 1, \ \varepsilon_j = -1$; 
\item[] $n_j \in N^{\circ} \overset{\mathrm{def}}{:=} \mathrm{Ker}E_D^N$   
when $\varepsilon_{j-1}  = -1, \ \varepsilon_j = 1$.  
\end{itemize}
We write the triple in the following way: 
\begin{equation*} 
\big(M,E_N^M,u(\theta)\big) := \big(N,E_D^N\big)\underset{D}{\bigstar}\big(\theta, E_{\theta(D)}^N\big). 
\end{equation*}  

\subsection{Modular Theory} Let $\psi$ be a faithful normal semifinite weight on $D$. Then the modular automorphism $\sigma_t^{\psi\circ E_D^N\circ E_N^M}$, $t \in \mathbf{R}$, is completely determined by 
\begin{equation}\label{eq1}
\sigma_t^{\psi\circ E_D^N\circ E_N^M}\left(u(\theta)\right)
= u(\theta) \left[D\psi\circ\theta^{-1}\circ E_{\theta(D)}^N : D\psi\circ E_D^N\right]_t. 
\end{equation}
See \cite[Theorem 4.1]{Ueda:JFA05}. In particular, if $N$ has a faithful normal trace $\tau$ satisfying that (i) both $E_D^N$ and $E_{\theta(D)}^N$ are $\tau$-preserving and (ii) $\tau|_{\theta(D)}\circ\theta = \tau|_D$, then the new state or weight $\tau\circ E_N^M$ ($= \tau|_D\circ E_D^N\circ E_N^M = \tau|_{\theta(D)}\circ E_{\theta(D)}^N\circ E_N^M$) becomes again a trace on $M$. Let $\widetilde{M} = M\rtimes_{\sigma^{\psi\circ E_D^N\circ E_{\theta(D)}^M}}\mathbf{R} \supseteq \widetilde{N} = N\rtimes_{\sigma^{\psi\circ E_D^N}}\mathbf{R} \supseteq \widetilde{D} = D\rtimes_{\sigma^{\psi}}\mathbf{R}$ be the inclusions of (continuous) {\it cores} with common canonical generators $\lambda(t)$, $t \in \mathbf{R}$, and then the {\it canonical liftings} $\widehat{E}_N^M : \widetilde{M} \rightarrow \widetilde{N}$, $\widehat{E}_D^N : \widetilde{N} \rightarrow \widetilde{D}$ of $E_N^M$, $E_D^N$ are provided in such a way that $\widehat{E}_N^M\big|_M=E_N^M$ and $\widehat{E}_D^N\big|_N=E_D^N$. Also, let $\widetilde{\theta} : \widetilde{D} \rightarrow \widetilde{N}$ be the {\it canonical extension} of $\theta$ defined by $\widetilde{\theta}|_D = \theta$ and $\widetilde{\theta}(\lambda(t)) = \left[\psi\circ\theta^{-1}\circ E_{\theta(D)}^N : \psi\circ E_D^N\right]_t \lambda(t)$ for $t \in \mathbf{R}$, and hence $\widetilde{\theta}\big(\widetilde{D}\big) = \widetilde{\theta(D)} = \theta(D)\rtimes_{\sigma^{\psi\circ\theta^{-1}}}\mathbf{R}$ so that we have the canonical lifting $\widehat{E}_{\theta(D)}^{N} : \widetilde{N} \rightarrow \widetilde{\theta}\big(\widetilde{D}\big)$ as before. Then, $\big(\widetilde{M},\widehat{E}_N^M,u(\theta)\big)$ is naturally identified with the HNN extension $\big(\widetilde{N},\widehat{E}_D^N\big)\bigstar_{\widetilde{D}}\big(\tilde{\theta},\widehat{E}_{\theta(D)}^N\big)$. See \cite[\S4]{Ueda:JFA05} for details. 

\subsection{$C^*$-Algebra Setup} Let $B \supseteq C$ be a unital inclusion of $C^*$-algebras, $\theta : C \rightarrow B$ be an injective unital $*$-homomorphism, and $E_C^B : B \rightarrow C$, $E_{\theta(C)}^B : B \rightarrow \theta(C)$ be conditional expectations. Assume, as a natural (or usual) requirement, that $E_C^B$ and $E_{\theta(C)}^B$ are non-degenerate (or equivalently have the faithful GNS representations), which ensures that $B$ is embedded in the reduced HNN extension faithfully. The {\it reduced HNN extension} $\big(B, E_C^B\big)\bigstar_C\big(\theta, E_{\theta(C)}^B\big)$ is constructed and defined as a triple $\big(A,E_B^A : A \rightarrow B,u(\theta)\big)$ in the exactly same manner as in the von Neumann algebra case, and it is indeed characterized by the same conditions (A), (M)  under the additional assumption that $E_B^A$ are non-degenerate. See \cite[\S\S7.2]{Ueda:JFA05} for the details. (An important thing about the characterization (\cite[Proposition 7.1]{Ueda:JFA05}) will be discussed in Remark \ref{rem2.3}.) In the $C^*$-algebra setup, another kind of HNN extension is available, and it is the {\it universal HNN extension} $B\bigstar_C^{\mathrm{univ}}\theta$, i.e., the universal $C^*$-algebra generated by $B$ and a single unitary $u(\theta)$ with subject to only the algebraic relations $u(\theta)\theta(c)u(\theta)^*=c$ for all $c \in C$.     

\section{Observation} 

\subsection{von Neumann Algebra Setup} Let $\big(M,E_N^M,u(\theta)\big) = \big(N,E_D^N\big) \bigstar_{D} \big(\theta, E_{\theta(D)}^N\big)$ be an HNN extension of von Neumann algebras, and 
\begin{equation*} 
\big(\mathcal{M},\mathcal{E}\big) := \big(N\otimes M_2(\mathbf{C}),E_{\theta} : \iota_{\theta}\big) \underset{D\oplus D }{\bigstar}\big(D\otimes  M_2(\mathbf{C}),E_1 : \iota_1\big)
\end{equation*} 
be the amalgamated free product of von Neumann algebra with the canonical embedding maps $\lambda, \lambda_{\theta}, \lambda_1$ of $D\oplus D, N\otimes M_2(\mathbf{C}), D\otimes M_2(\mathbf{C})$, respectively, into $\mathcal{M}$, where 
\begin{gather*} 
\iota_{\theta}\left((d_1, d_2)\right) := 
\begin{bmatrix} d_1 & 0 \\ 0 & \theta(d_2) \end{bmatrix}, \quad 
\iota_1\left((d_1, d_2)\right) := 
\begin{bmatrix} d_1 & 0 \\ 0 & d_2 \end{bmatrix}; \\ 
E_{\theta} := \begin{bmatrix} E_D^N & 0 \\ 0 & E_{\theta(D)}^N \end{bmatrix}, \quad 
E_1 := \mathrm{Id}\otimes E_{\mathbf{C}^2}^{M_2(\mathbf{C})}; \\
\lambda = \lambda_{\theta}\circ\iota_{\theta} = \lambda_1\circ\iota_1. 
\end{gather*} 
See \cite[\S2]{Ueda:JFA05} for the construction and terminologies. Let us denote by $\mathcal{E}_{\theta}$ the conditional expectation from $\mathcal{M}$ onto $\lambda_{\theta}\left(N\otimes M_2(\mathbf{C})\right)$ that satisfies $\mathcal{E}\circ\mathcal{E}_{\theta} = \mathcal{E}$. 

\begin{proposition}\label{prop2.1}
There is a bijective $*$-homomorphism $\Phi : \mathcal{M} \rightarrow M\otimes M_2(\mathbf{C})$ such that $\Phi\left(\lambda_{\theta}\left(N\otimes M_2(\mathbf{C})\right)\right) = N\otimes  M_2(\mathbf{C}) \subseteq M\otimes M_2(\mathbf{C})$ and moreover 
\begin{equation}\label{eq2}
\Phi\circ\mathcal{E}_{\theta} = \left(E_N^M\otimes\mathrm{Id}\right)\circ\Phi.  
\end{equation}
The above bijective $*$-homomorphism $\Phi$ is precisely given by 
\begin{equation*}
\Phi : \left\{
\begin{array}{ccc} 
\lambda_1\left(\begin{bmatrix} 0 & 1 \\ 0 & 0 \end{bmatrix}\right)\lambda_{\theta}\left(\begin{bmatrix} 0 & 0 \\ 1 & 0 \end{bmatrix}\right) & \longmapsto & \begin{bmatrix} u(\theta) & 0 \\ 0 & 0 \end{bmatrix},  \\ \\ 
\lambda_{\theta}\left(\begin{bmatrix} n & 0 \\ 0 & 0 \end{bmatrix}\right) & \longmapsto & \begin{bmatrix} n & 0 \\ 0 & 0 \end{bmatrix}, \\ \\ 
\lambda_{\theta}\left(\begin{bmatrix} 0 & 1 \\ 0 & 0 \end{bmatrix}\right) & \longmapsto & \begin{bmatrix} 0 & 1 \\ 0 & 0 \end{bmatrix}.
\end{array}
\right.
\end{equation*} 
\end{proposition}  
\begin{proof} 
Let us first recall (and improve) the construction of reduced HNN extensions given in \cite{Ueda:JFA05}. Let $(\mathcal{M},\mathcal{E})$ be as in the statement, and the HNN extension $\big(M,E_N^M,u(\theta)\big)$ is realized in the compressed algebra $p\mathcal{M}p$ with $p := \lambda_{\theta}\left(\begin{bmatrix} 1 & 0 \\ 0 & 0 \end{bmatrix}\right)$ as follows. (Note that another algebra slightly larger than $\mathcal{M}$ was used in \cite{Ueda:JFA05}, but it is clear that  $\mathcal{M}$ is sufficiently large to construct the desired algebra.) Identify $n \in N$ with $\lambda_{\theta}\left(\begin{bmatrix} n & 0 \\ 0 & 0 \end{bmatrix}\right)$
and set $u(\theta) := \lambda_1\left(\begin{bmatrix} 0 & 1 \\ 0 & 0 \end{bmatrix}
\right)\lambda_{\theta}\left(\begin{bmatrix} 0 & 0 \\ 1 & 0 \end{bmatrix}\right)$, and then the desired algebra $M$ is generated by $N$ and $u(\theta)$ inside $p\mathcal{M}p$, and the conditional expectation $E_N^M$ is obtained as the restriction of $\mathcal{E}_{\theta}$ to $M$. 

Let $\Phi : \mathcal{M} \rightarrow p\mathcal{M}p\otimes M_2(\mathbf{C})$ be the bijective normal $*$-isomorphism determined by the $2\times2$ matrix unit system 
\begin{equation*} 
p= \lambda_{\theta}\left(\begin{bmatrix} 1 & 0 \\ 0 & 0 \end{bmatrix}\right),\quad 
\lambda_{\theta}\left(\begin{bmatrix} 0 & 1 \\ 0 & 0 \end{bmatrix}\right),\quad  
\lambda_{\theta}\left(\begin{bmatrix} 0 & 0 \\ 1 & 0 \end{bmatrix}\right),\quad 
\lambda_{\theta}\left(\begin{bmatrix} 0 & 0 \\ 0 & 1 \end{bmatrix}\right). 
\end{equation*} 
Then, we get 
\begin{equation*} 
\Phi\left(u(\theta)\right) = \begin{bmatrix} u(\theta) & 0 \\ 0 & 0 \end{bmatrix};\quad 
\Phi\left(\lambda_{\theta}\left(\begin{bmatrix} n & 0 \\ 0 & 0 \end{bmatrix}\right)\right) = \begin{bmatrix} n & 0 \\ 0 & 0 \end{bmatrix};\quad 
\Phi\left(\lambda_{\theta}\left(\begin{bmatrix} 0 & 1 \\ 0 & 0 \end{bmatrix}\right)\right) = \begin{bmatrix} 0 & 1 \\ 0 & 0 \end{bmatrix}, 
\end{equation*} 
where the right-hand sides are considered in $M_2(p\mathcal{M}p) = p\mathcal{M}p\otimes M_2(\mathbf{C})$. Note that 
\begin{equation*} 
u(\theta) = \lambda_1\left(\begin{bmatrix} 0 & 1 \\ 0 & 0 \end{bmatrix}\right)\lambda_{\theta}\left(\begin{bmatrix} 0 & 0 \\ 1 & 0 \end{bmatrix}\right),\quad 
\lambda_{\theta}\left(\begin{bmatrix} n & 0 \\ 0 & 0 \end{bmatrix}\right), \quad \lambda_{\theta}\left(\begin{bmatrix} 0 & 1 \\ 0 & 0 \end{bmatrix}\right)
\end{equation*} 
(considered in $\mathcal{M}$) generate the whole $\mathcal{M}$ since 
\begin{equation*} 
\lambda_1\left(\begin{bmatrix} 0 & 1 \\ 0 & 0 \end{bmatrix}\right)\lambda_{\theta}\left(\begin{bmatrix} 0 & 0 \\ 1 & 0 \end{bmatrix}\right) \times \lambda_{\theta}\left(\begin{bmatrix} 0 & 1 \\ 0 & 0 \end{bmatrix}\right) = \lambda_1\left(\begin{bmatrix} 0 & 1 \\ 0 & 0 \end{bmatrix}\right). 
\end{equation*} 
Since  $\Phi(\mathcal{M}) = p\mathcal{M}p\otimes M_2(\mathbf{C})$, we conclude that $p\mathcal{M}p\otimes M_2(\mathbf{C}) = M\otimes M_2(\mathbf{C})$ and $M = p\mathcal{M}p$. The equality \eqref{eq2} is easily verified. 
\end{proof}

Remark that the above argument clearly works well even in the general case where the $\theta$ is replaced by a family of injective normal unital $*$-homomorphisms from $D$  into $N$. 

\subsection{Reduced $C^*$-Algebra Setup} Let $\big(A,E_B^A,u(\theta)\big) = \big(B,E_C^B\big) \bigstar_{C} \big(\theta, E_{\theta(C)}^B\big)$ be a reduced HNN extension of $C^*$-algebras, and 
\begin{equation*} 
\left(\mathcal{A},\mathcal{E}\right) := \left(B\otimes M_2(\mathbf{C}),E_{\theta} : \iota_{\theta}\right) \underset{C\oplus C}{\bigstar} \left(C\otimes M_2(\mathbf{C}),E_1 : \iota_1\right)
\end{equation*} 
be the reduced amalgamated free product of $C^*$-algebras with the canonical embedding maps $\lambda, \lambda_{\theta}, \lambda_1$ of $C\oplus C, B\otimes M_2(\mathbf{C}), C\otimes M_2(\mathbf{C})$, respectively, into $\mathcal{M}$, where 
\begin{gather*} 
\iota_{\theta}\left((c_1, c_2)\right) := 
\begin{bmatrix} c_1 & 0 \\ 0 & \theta(c_2) \end{bmatrix}, \quad 
\iota_1\left((c_1, c_2)\right) := 
\begin{bmatrix} c_1 & 0 \\ 0 & c_2 \end{bmatrix}; \\ 
E_{\theta} := \begin{bmatrix} E_C^B & 0 \\ 0 & E_{\theta(C)}^B \end{bmatrix}, \quad 
E_1 := \mathrm{Id}\otimes E_{\mathbf{C}^2}^{M_2(\mathbf{C})}; \\
\lambda = \lambda_{\theta}\circ\iota_{\theta} = \lambda_1\circ\iota_1. 
\end{gather*} 
Denote by $\mathcal{E}_{\theta}$ the conditional expectation from $\mathcal{A}$ onto $\lambda_{\theta}\left(B\otimes M_2(\mathbf{C})\right)$ that satisfies $\mathcal{E}\circ\mathcal{E}_{\theta} = \mathcal{E}$. (See e.g.~\cite[Lemma 1.1]{BlanchardDykema:Pacific01}.) The proposition below is shown in the exactly same way as in the von Neumann algebra setting. 

\begin{proposition}\label{prop2.2}
There is a bijective $*$-homomorphism $\Phi : \mathcal{A} \rightarrow A\otimes M_2(\mathbf{C})$ such that $\Phi\left(\lambda_{\theta}\left(B\otimes M_2(\mathbf{C})\right)\right) = B\otimes  M_2(\mathbf{C}) \subseteq A\otimes M_2(\mathbf{C})$ and moreover 
\begin{equation*}
\Phi\circ\mathcal{E}_{\theta} = \left(E_B^A\otimes\mathrm{Id}\right)\circ\Phi.  
\end{equation*}
The above bijective $*$-homomorphism $\Phi$ is precisely given by 
\begin{equation*}
\Phi : \left\{
\begin{array}{ccc} 
\lambda_1\left(\begin{bmatrix} 0 & 1 \\ 0 & 0 \end{bmatrix}\right)\lambda_{\theta}\left(\begin{bmatrix} 0 & 0 \\ 1 & 0 \end{bmatrix}\right) & \longmapsto & \begin{bmatrix} u(\theta) & 0 \\ 0 & 0 \end{bmatrix},  \\ \\ 
\lambda_{\theta}\left(\begin{bmatrix} b & 0 \\ 0 & 0 \end{bmatrix}\right) & \longmapsto & \begin{bmatrix} b & 0 \\ 0 & 0 \end{bmatrix}, \\ \\ 
\lambda_{\theta}\left(\begin{bmatrix} 0 & 1 \\ 0 & 0 \end{bmatrix}\right) & \longmapsto & \begin{bmatrix} 0 & 1 \\ 0 & 0 \end{bmatrix}.
\end{array}
\right.
\end{equation*} 
%In particular, the HNN extension $A=B\bigstar_{C}\theta$ is isomorphic to the corner algebra $p\mathcal{A}p$ of the amalgamated free product $\mathcal{A}=\left(B\otimes M_2(\mathbf{C})\right)\bigstar_{C\otimes\mathbf{C}^2}\left(C\otimes M_2(\mathbf{C})\right)$ with $p := \lambda_{\theta}\left(\begin{bmatrix} 1 & 0 \\ 0 & 0 \end{bmatrix}\right)$, and moreover $\mathcal{A}$ is isomorphic to $A\otimes M_2(\mathbf{C})$.  
\end{proposition}

The statement of Proposition \ref{prop2.2} still holds true even when the $\theta$ is replaced by a family $\Theta$ of injective unital $*$-homomorphisms from $C$ into $B$,  but when $\Theta$ is an infinite family one has to replace $\mathcal{A}$ by the $C^*$-subalgebra generated by $\lambda_{\Theta}(B\otimes\mathbb{K}(\ell^2(\Theta_1)))$ and $\lambda_1(C\otimes\mathbb{K}(\ell^2(\Theta_1)))$ with the notations in \cite[\S7]{Ueda:JFA05}, where $\mathbb{K}(\mathcal{H})$ denotes the algebra of all compact operators on a Hilbert space $\mathcal{H}$. The proof of Proposition \ref{prop2.1} still works without change when $\Theta$ is a finite family. The case when $\Theta$ is infinite needs to pass through the inductive limit by finite subfamilies $\Xi \nearrow \Theta$ with the aid of \cite[Theorem 1.3]{BlanchardDykema:Pacific01}.  

\begin{remark}\label{rem2.3}{\rm There is an insufficient point related to the characterization of reduced HNN extensions (\cite[Proposition 7.1]{Ueda:JFA05}); in fact, we did not prove that the reduced HNN extensions that we constructed in \cite{Ueda:JFA05} actually satisfy the condition (ii) (the non-degeneracy condition) there. Of course, this is not a problem in several cases including reduced group $C^*$-algebras associated with HNN extensions of groups. However, it is certainly necessary to prove it for the justification of our definition. One easy way to do so is provided by Proposition \ref{prop2.2} as follows. Notice that the proof of Proposition \ref{prop2.2} (or Proposition \ref{prop2.1}) was done based on the construction. Thus, Proposition \ref{prop2.2} shows that $E_B^A\otimes\mathrm{Id} : A\otimes M_2(\mathbf{C}) \rightarrow B\otimes M_2(\mathbf{C})$ is non-degenerate (since so is $\mathcal{E}_{\theta}$ by the amalgamated free product construction), which immediately implies that so is $E_B^A$. Note that we used in \cite{Ueda:JFA05} a reduced amalgamated free product larger than the above $\mathcal{A}$ to construct the reduced HNN extension $A$, and thus it is necessary to prove that this $A$ is the same as that constructed there without the use of \cite[Proposition 7.1]{Ueda:JFA05}. However, this is not a problem because $\mathcal{A}$ is naturally embedded into the larger one faithfully thanks to \cite[Theorem 1.3]{BlanchardDykema:Pacific01}. \\
} 
\end{remark}     

\subsection{Universal $C^*$-Algebra Setup} Let $B \supseteq C$ be a unital inclusion of $C^*$-algebras with an injective unital $*$-homomorphism $\theta : C \rightarrow B$ as above, and $A = B\bigstar_C^{\mathrm{univ}} \theta$ be the universal HNN extension of $C^*$-algebras. Let 
\begin{equation*} 
\mathcal{A} := \left(B\otimes M_2(\mathbf{C}) : \iota_{\theta}\right)\underset{C\oplus C}{\bigstar^{\mathrm{univ}}} \left(C\otimes M_2(\mathbf{C}) : \iota_1\right)
\end{equation*} 
be the universal amalgamated free product of $C^*$-algebras over $C\oplus C$ via the distinguished embedding maps 
\begin{equation*} 
\iota_{\theta}\left((c_1, c_2)\right) := \begin{bmatrix} c_1 & 0 \\ 0 & \theta(c_2) \end{bmatrix}, \quad \iota_1\left((c_1, c_2)\right) := \begin{bmatrix} c_1 & 0 \\ 0 & c_2 \end{bmatrix}. 
\end{equation*} 
Let us denote by $j$, $j_{\theta}$ and $j_1$ the canonical embedding maps of $C\oplus C$, $B\otimes M_2(\mathbf{C})$ and $C\otimes M_2(\mathbf{C})$ into $\mathcal{A}$, respectively, that satisfy $j=j_{\theta}\circ\iota_{\theta} = j_1\circ\iota_1$.  

\begin{proposition}\label{prop2.4} There is a bijective $*$-homomorphism $\Phi : \mathcal{A} \rightarrow A\otimes M_2(\mathbf{C})$ such that 
\begin{equation}
\Phi : \left\{
\begin{array}{ccc} 
j_1\left(\begin{bmatrix} 0 & 1 \\ 0 & 0 \end{bmatrix}\right)j_{\theta}\left(\begin{bmatrix} 0 & 0 \\ 1 & 0 \end{bmatrix}\right) & \longmapsto & \begin{bmatrix} u(\theta) & 0 \\ 0 & 0 \end{bmatrix},  \\ \\ 
j_{\theta}\left(\begin{bmatrix} b & 0 \\ 0 & 0 \end{bmatrix}\right) & \longmapsto & \begin{bmatrix} b & 0 \\ 0 & 0 \end{bmatrix}, \\ \\ 
j_{\theta}\left(\begin{bmatrix} 0 & 1 \\ 0 & 0 \end{bmatrix}\right) & \longmapsto & \begin{bmatrix} 0 & 1 \\ 0 & 0 \end{bmatrix}.
\end{array}
\right.\label{eq3}
\end{equation} 
\end{proposition}     
\begin{proof} 
Let us first define two $*$-homomorphisms $\Phi_{\theta} : B\otimes M_2(\mathbf{C}) \rightarrow A\otimes M_2(\mathbf{C})$, $\Phi_1 : C\otimes M_2(\mathbf{C}) \rightarrow A\otimes M_2(\mathbf{C})$ by 
\begin{align*} 
\Phi_{\theta}\left(j_{\theta}\left(\begin{bmatrix} b_{11} & b_{12} \\ b_{21} & b_{22} \end{bmatrix}\right)\right) &:= 
\begin{bmatrix} b_{11} & b_{12} \\ b_{21} & b_{22} \end{bmatrix}; \\
\Phi_1\left(j_1\left(\begin{bmatrix} c_{11} & c_{12} \\ c_{21} & c_{22} \end{bmatrix}\right)\right) &:= \begin{bmatrix} c_{11} & c_{12} u(\theta) \\ u(\theta)^* c_{21} & \theta\left(c_{22}\right) \end{bmatrix}\, \left(= \begin{bmatrix} 1 & 0 \\ 0 & u(\theta)^* \end{bmatrix}\begin{bmatrix} c_{11} & c_{12} \\ c_{21} & c_{22} \end{bmatrix}\begin{bmatrix} 1 & 0 \\ 0 & u(\theta) \end{bmatrix}\right). 
\end{align*} 
Then, we have 
\begin{gather*} 
\Phi_{\theta}\left(j_{\theta}\circ \iota_{\theta}\left((c_1, c_2)\right)\right) = \Phi_{\theta}\left(j_{\theta}\left(\begin{bmatrix} c_{11} & 0 \\ 0 & \theta\left(c_{22}\right) \end{bmatrix}\right)\right) = \begin{bmatrix} c_{11} & 0 \\ 0 & \theta\left(c_{22}\right) \end{bmatrix}, \\
\Phi_1\left( j_1\circ\iota_1\left((c_1, c_2)\right)\right) = \Phi_1\left(j_1\left(\begin{bmatrix} c_{11} & 0 \\ 0 & c_{22} \end{bmatrix}\right)\right) = \begin{bmatrix} c_{11} & 0 \\ 0 & \theta\left(c_{22}\right) \end{bmatrix}.  
\end{gather*}  
Thus the universality of $\displaystyle{\mathcal{A} = \left(B\otimes M_2(\mathbf{C}) : \iota_{\theta}\right)\underset{C\oplus C}{\bigstar} \left(C\otimes M_2(\mathbf{C}) : \iota_1\right)}$ ensures that there is a unique unital $*$-homomorphism $\Phi := \Phi_{\theta}\bigstar\Phi_1 : \mathcal{A} \rightarrow A\otimes M_2(\mathbf{C})$ extending both $\Phi_{\theta}$ and $\Phi$. Since $\Phi$ agrees with \eqref{eq3}, it remains only to show that $\Phi$ is bijective. To do so we will construct the inverse of $\Phi$ in what follows. By the universality of $A = B\bigstar_C^{\mathrm{univ}} \theta$ we can construct the unital $*$-homomorphism $\Psi_0 : A \rightarrow p\mathcal{A}p$ with $p := j_{\theta}\left(\begin{bmatrix} 1 & 0 \\ 0 & 0 \end{bmatrix}\right)$ in such a way that 
\begin{equation*} 
\Psi_0(b) := j_{\theta}\left(\begin{bmatrix} b & 0 \\ 0 & 0 \end{bmatrix}\right), \quad 
\Psi_0\left(u(\theta)\right) := j_1\left(\begin{bmatrix} 0 & 1 \\ 0 & 0 \end{bmatrix}\right)j_{\theta}\left(\begin{bmatrix} 0 & 0 \\ 1 & 0 \end{bmatrix}\right)
\end{equation*}  
since 
\begin{equation*} 
j_1\left(\begin{bmatrix} 0 & 1 \\ 0 & 0 \end{bmatrix}\right)j_{\theta}\left(\begin{bmatrix} 0 & 0 \\ 1 & 0 \end{bmatrix}\right)j_{\theta}\left(\begin{bmatrix} \theta(c) & 0 \\ 0 & 0 \end{bmatrix}\right)j_{\theta}\left(\begin{bmatrix} 0 & 1 \\ 0 & 0 \end{bmatrix}\right)j_1\left(\begin{bmatrix} 0 & 0 \\ 1 & 0 \end{bmatrix}\right) = j_{\theta}\left(\begin{bmatrix} c & 0 \\ 0 & 0 \end{bmatrix}\right) 
\end{equation*} 
for all $c \in C$. Consider the following two $2\times2$ matrix unit systems inside $\mathbf{C}1\otimes M_2(\mathbf{C})\subseteq A\otimes M_2(\mathbf{C})$ and $\mathcal{A}$ 
\begin{gather*} 
e_{11} := \begin{bmatrix} 1 & 0 \\ 0 & 0 \end{bmatrix},\ 
e_{12} := \begin{bmatrix} 0 & 1 \\ 0 & 0 \end{bmatrix},\ 
e_{21} := \begin{bmatrix} 0 & 0 \\ 1 & 0 \end{bmatrix},\ 
e_{22} := \begin{bmatrix} 0 & 0 \\ 0 & 1 \end{bmatrix}; \\
f_{11} := j_{\theta}\left(\begin{bmatrix} 1 & 0 \\ 0 & 0 \end{bmatrix}\right),\ 
f_{12} := j_{\theta}\left(\begin{bmatrix} 0 & 1 \\ 0 & 0 \end{bmatrix}\right),\ 
f_{21} := j_{\theta}\left(\begin{bmatrix} 0 & 0 \\ 1 & 0 \end{bmatrix}\right),\ 
f_{22} := j_{\theta}\left(\begin{bmatrix} 0 & 0 \\ 0 & 1 \end{bmatrix}\right), 
\end{gather*}
respectively, with $f_{11} = p$, and then $\Psi_0$ is extended to a $*$-homomorphism $\Psi: A\otimes M_2(\mathbf{C}) \rightarrow \mathcal{A}$ by $\Psi(x) := 
\sum_{i,j=1}^2 f_{i1}\Psi_0(e_{1i}x e_{j1})f_{1j}$  
%\begin{align*}
%\Psi(x) := 
%f_{11}\Psi_0\left(e_{11}x e_{11}\right)f_{11} &+  f_{11}\Psi_0\left(e_{11} x e_{21}\right)f_{12} \\ 
%&+ f_{21}\Psi_0\left(e_{12} x e_{11}\right)f_{11} +  f_{21}\Psi_0\left(e_{12} x e_{21}\right)f_{12}
%\end{align*} 
for $x \in A\otimes M_2(\mathbf{C})$. Then one immediately observes  
\begin{align*} 
\Psi\left(\begin{bmatrix} b & 0 \\ 0 & 0 \end{bmatrix}\right) 
%&= f_{11}\Psi_0\left(\begin{bmatrix} b & 0 \\ 0 & 0 \end{bmatrix}\right)f_{11} \\
%&= j_{\theta}\left(\begin{bmatrix} 0 & 0 \\ 0 & 1 \end{bmatrix}\right)j_{\theta}\left(\begin{bmatrix} 0 & 0 \\ 0 & b \end{bmatrix}\right) j_{\theta}\left(\begin{bmatrix} 0 & 0 \\ 0 & 1 \end{bmatrix}\right) \\
&= j_{\theta}\left(\begin{bmatrix} b & 0 \\ 0 & 0 \end{bmatrix}\right), \\
\Psi\left(\begin{bmatrix} u(\theta) & 0 \\ 0 & 0 \end{bmatrix}\right) 
%&= f_{11}\Psi_0\left(\begin{bmatrix} u(\theta) & 0 \\ 0 & 0 \end{bmatrix}\right)f_{11} \\
%&= j_{\theta}\left(\begin{bmatrix} 0 & 0 \\ 0 & 1 \end{bmatrix}\right)
%j_{\theta}\left(\begin{bmatrix} 0 & 0 \\ 1 & 0 \end{bmatrix}\right)j_1\left(\begin{bmatrix} 0 & 1 \\ 0 & 0 \end{bmatrix}\right)j_{\theta}\left(\begin{bmatrix} 0 & 0 \\ 0 & 1 \end{bmatrix}\right) \\
%&= j_{\theta}\left(\begin{bmatrix} 0 & 0 \\ 0 & 1 \end{bmatrix}\right)
%j_{\theta}\left(\begin{bmatrix} 0 & 0 \\ 1 & 0 \end{bmatrix}\right)j_1\left(\begin{bmatrix} 0 & 1 \\ 0 & 0 \end{bmatrix}\right)j_1\left(\begin{bmatrix} 0 & 0 \\ 0 & 1 \end{bmatrix}\right) \\
&=j_1\left(\begin{bmatrix} 0 & 1 \\ 0 & 0 \end{bmatrix}\right) j_{\theta}\left(\begin{bmatrix} 0 & 0 \\ 1 & 0 \end{bmatrix}\right), \\
\Psi\left(\begin{bmatrix} 0 & 1 \\ 0 & 0 \end{bmatrix}\right) 
%&= \Psi\left(e_{12}\right) \\
%&= f_{11} \Psi_0\left(e_{11}e_{12}e_{21}\right)f_{12} \\
%&= f_{11} f_{11} f_{12} \\ 
&= f_{12}\, \left(= j_{\theta}\left(\begin{bmatrix} 0 & 1 \\ 0 & 0 \end{bmatrix}\right)\right),  
\end{align*}   
and hence  $\Psi\circ\Phi = \mathrm{id}_{\mathcal{A}}$ and $\Phi\circ\Psi = \mathrm{id}_{A\otimes M_2(\mathbf{C})}$.
\end{proof}

The statement of Proposition \ref{prop2.4} still holds true even when the $\theta$ is replaced by a family $\Theta$ of injective unital $*$-homomorphisms from $C$ into $B$, but the same care as in the reduced setting is required. Also, it should be pointed out that (the first half of) the above proof says that the matrix trick we employ provides a simple way to construct universal HNN extensions of $C^*$-algebras.

\subsection{Equivalence Relation Case: Relation to Gaboriau's work}

In \cite{Gaboriau:InventMath00} (see also a related and a bit earlier work due to Paulin \cite{Paulin:MarkovProcess99}) Gaboriau introduced the notion of HNN equivalence relations and derive a formula of costs for them from that for amalgamated free product equivalence relations based on a certain relationship between HNN and amalgamated free product equivalence relations, which we will explain below. 

Let $\mathcal{R} \supseteq \mathcal{S} \supseteq \mathcal{T}$ be discrete standard Borel equivalence relations over a standard Borel space $X$, and $\varTheta : E_{+1} \rightarrow E_{-1}$ be a Borel isomorphism, called a {\it partial transformation}, between two Borel subsets $E_{+1}, E_{-1} \subseteq X$. Assume that $\mathcal{T}$ is trivial on $X\setminus E_{-1}$, i.e., if $x \in X\setminus E_{-1}$ then $y=x$ holds for every $(x,y) \in \mathcal{T}$. Set $\mathcal{T}_{\varTheta} := \left\{ (x,x) : x \in X\setminus E_{+1} \right\} \sqcup \left\{ (x,y) \in E_{+1}\times E_{+1} : (\varTheta(x),\varTheta(y)) \in \mathcal{T} \right\}$, a new equivalence relation over $X$, and suppose $\mathcal{T}_{\varTheta} \subseteq \mathcal{S}$. We write $x \overset{\varTheta^{\varepsilon}}{\rightarrow} y$ for $x,y \in X$ and $\varepsilon = \pm1$ when $x \in E_{\varepsilon}$ and $y=\varTheta^{\varepsilon}(x)$ or equivalently $(x,y) \in \mathrm{Graph}(\varTheta^{\varepsilon}) := \left\{ (x,\varTheta^{\varepsilon}(x)) \in X\times X : x \in E_{\varepsilon} \right\}$ (with letting $\varTheta^{+1} := \varTheta$). A finite sequence $(x_1,\dots,x_{2n})$ in $X$ with 
\begin{equation}\label{eq4}
x_1 \underset{\mathcal{S}}{\sim}\quad \cdots x_{2i-2} \overset{\varTheta^{\varepsilon_{i-1}}}{\rightarrow} x_{2i-1} \underset{\mathcal{S}}{\sim} x_{2i} \overset{\varTheta^{\varepsilon_{i}}}{\rightarrow} x_{2i+1} \quad \cdots \underset{\mathcal{S}}{\sim} x_{2n} 
\end{equation} 
is called a reduced word if $n \geq 2$ and no subsequence $(x_{2i-2},x_{2i-1},x_{2i})$ satisfying  
\begin{align*} 
x_{2i-2} \overset{\varTheta}{\rightarrow} x_{2i-1} \underset{\mathcal{S}}{\sim} x_{2i} \overset{\varTheta^{-1}}{\rightarrow} x_{2i+1} &\quad \text{and} \quad (x_{2i-1},x_{2i}) \in \mathcal{T} \quad  \text{nor} \\
x_{2i-2} \overset{\varTheta^{-1}}{\rightarrow} x_{2i-1} \underset{\mathcal{S}}{\sim} x_{2i} \overset{\varTheta}{\rightarrow} x_{2i+1} &\quad \text{and} \quad (x_{2i-1},x_{2i}) \in \mathcal{T}_{\varTheta}
\end{align*} 
appears in \eqref{eq4}; or if $n=1$ and $x_1 \neq x_2$. Then, the bigger $\mathcal{R}$ is said to be the {\it HNN extension of $\mathcal{S}$ by $\varTheta$} and denoted by $\mathcal{S}\bigstar_{\mathcal{T}}\varTheta$, if  
\begin{itemize}
\item $\mathcal{R}$ is generated by $\mathcal{S}$ and $\varTheta$ or more precisely the smallest equivalence relation containing $\mathcal{S}\cup\mathrm{Graph}(\varTheta)$;  
\item any reduced word $(x_1,\dots,x_{2n})$ with $x_1,\dots,x_{2n} \in X$ (in the above sense) must satisfy that $x_1 \neq x_{2n}$. 
\end{itemize}
In the measurable setting, i.e., $X$ is equipped with a regular Borel measure $\mu$, all the equivalence relations and $\varTheta$ are assumed to be non-singular under $\mu$, and $\mathcal{R} = \mathcal{S}\bigstar_{\mathcal{T}}\varTheta$ is defined in the same way but up to $\mu$-null set, i.e., there is a $\mu$-conull subset of $X$, on which the condition holds. In what follows, we consider in the measurable setting and only the case that $E_{\pm1} = X$; namely $\varTheta$ is a transformation defined on the entire space $X$. In this case, simply $\mathcal{T}_{\varTheta} = (\varTheta\times\varTheta)^{-1}(\mathcal{T})$. Let $W^*(\mathcal{R}) \supseteq W^*(\mathcal{S}) \supseteq W^*(\mathcal{T}),\, W^*(\mathcal{T}_{\varTheta}) \left(\supseteq L^{\infty}(X)\right)$ be the von Neumann algebras associated with $\mathcal{R} \supseteq \mathcal{S} \supseteq \mathcal{T},\, \mathcal{T}_{\varTheta}$ constructed by Feldman-Moore's construction \cite{FeldmanMoore:TAMS77-2}. Let $E_{\mathcal{S}}^{\mathcal{R}} : W^*(\mathcal{R}) \rightarrow W^*(\mathcal{S})$, $E_{\mathcal{T}}^{\mathcal{S}} : W^*(\mathcal{S}) \rightarrow W^*(\mathcal{T})$, $E_{\mathcal{T}_{\varTheta}}^{\mathcal{S}} : W^*(\mathcal{S}) \rightarrow W^*(\mathcal{T}_{\varTheta})$ be the unique (faithful normal) conditional expectations. Let us also denote by $\lambda$ the left regular representation of $\mathcal{R}$ in the Feldman-Moore construction, which in particular gives a representation $\lambda(R)$ of each partial transformation $R$ with $\mathrm{Graph}(R) \subseteq \mathcal{R}$ on $L^2(\mathcal{R},\mu_r)$ with the right-counting measure $\mu_r$, by the convolution operator $\lambda\big(\chi_{\mathrm{Graph}(R^{-1})}\big)$ of the characteristic function $\chi_{\mathrm{Graph}(R)}(x,y)$ in the terminologies, see \cite[Proposition 2.1]{FeldmanMoore:TAMS77-2}. Note that the $\lambda(R)$'s generate the $W^*(\mathcal{R})$ as von Neumann algebra, and indeed this is an important one in the set of axioms of the notion of Cartan subalgberas (see \cite[Definition 3.1 and Theorem 1]{FeldmanMoore:TAMS77-2}). Here the following remarks are in order: (a) For $\mathcal{X} = \mathcal{S}, \mathcal{T}, \mathcal{T}_{\varTheta}$, $W^*(\mathcal{X})$ is the s.o.-closure of the linear span of all $\lambda(R)$ with $\mathrm{Graph}(R) \subseteq \mathcal{X}$ inside $W^*(\mathcal{R})$. (b) For each partial transformation $S$ with $\mathrm{Graph}(S) \subseteq \mathcal{S}$, one has $E^{\mathcal{S}}_{\mathcal{T}}\left(\lambda\big(\chi_{\mathrm{Graph}(S)}\big)\right) = \lambda\big(\chi_{\mathrm{Graph}(S)\cap\mathcal{X}}\big)$ for $\mathcal{X} = \mathcal{T}, \mathcal{T}_{\varTheta}$.  
%\begin{equation*}
%E^{\mathcal{S}}_{\mathcal{T}}\left(\lambda\big(\chi_{\mathrm{Graph}(S)}\big)\right) = \lambda\big(\chi_{\mathrm{Graph}(S)\cap\mathcal{X}}\big).   
%\end{equation*}

Here is an expected fact.     

\begin{proposition}\label{prop2.5} Set $u(\theta) := \lambda(\varTheta)$ and denote by $\theta$ the injective unital $*$-homomorphism from $W^*(\mathcal{T})$ onto $W^*(\mathcal{T}_{\varTheta})$ implemented by $u(\theta)^* = \lambda(\varTheta^{-1})$. Then, $\big(W^*(\mathcal{R}), E^{\mathcal{R}}_{\mathcal{S}}, u(\theta)\big)$ is identified with  
$\big(W^*(\mathcal{S}),E_{\mathcal{T}}^{\mathcal{S}}\big)\bigstar_{W^*(\mathcal{T})}\big(\theta, E_{\mathcal{T}_{\varTheta}}^{\mathcal{S}}\big)$.  
\end{proposition}  
\begin{proof} Since $\mathcal{R}$ is generated by $\mathcal{S}$ and $\mathrm{Graph}(\varTheta)$, it is plain to see that the $\lambda(S)$'s with $\mathrm{Graph}(S) \subseteq \mathcal{S}$ and $u(\theta)$ generate $W^*(\mathcal{R})$ as von Neumann algebra. Hence it suffices to confirm that the triple $\big(W^*(\mathcal{R}), E_{\mathcal{S}}^{\mathcal{R}}, u(\theta)\big)$ satisfies the condition (M) in \S\S1.1 because the condition (A) is trivial by the definition of $\theta$ and $u(\theta)$. Thanks to the above remarks (a),(b) what we need is to confirm that $E_{\mathcal{S}}^{\mathcal{R}}(w) = 0$ for any word $w = u(\theta)^{\varepsilon_0}\lambda(S_1)u(\theta)^{\varepsilon_1}\lambda(S_2)\cdots\lambda(S_{\ell})u(\theta)^{\varepsilon_{\ell}}$ with partial transformations $S_1,\dots,S_{\ell}$, satisfying that $\mathrm{Graph}(S_j) \subseteq \mathcal{S}$, $j=1,\dots,\ell$, and moreover that if $\varepsilon_{j-1} \neq \varepsilon_j$ implies that 
\begin{equation*} 
\begin{matrix} 
\mathrm{Graph}(S_j^{-1}) \subseteq \mathcal{S}\setminus\mathcal{T}_{\varTheta} & \text{when $\varepsilon_{j-1} = 1$, $\varepsilon_j = -1$}; \\
\mathrm{Graph}(S_j^{-1}) \subseteq \mathcal{S}\setminus\mathcal{T} & \text{when $\varepsilon_{j-1} = -1$, $\varepsilon_j = 1$}. 
\end{matrix}
\end{equation*} 
By the definition of convolution operators, it is easy to re-write the word $w$ to be 
\begin{equation*} 
\lambda\big(\varTheta^{\varepsilon_0}\circ S_1\circ\varTheta^{\varepsilon_1}\circ S_2 \circ\cdots\circ S_{\ell}\circ\varTheta^{\varepsilon_{\ell}}\big),  
\end{equation*} 
where $\varTheta^{\varepsilon_0}\circ S_1\circ\varTheta^{\varepsilon_1}\circ S_2 \circ\cdots\circ S_{\ell}\circ\varTheta^{\varepsilon_{\ell}}$ means the successive composition of partial transformations. 
Let $X_0$ be a $\mu$-conull subset of $X$, on which the condition of being $\mathcal{R} = \mathcal{S}\bigstar_{\mathcal{T}}\varTheta$ holds, and choose $(x,y) \in X_0\times X_0$ from the graph of $(\varTheta^{\varepsilon_0}\circ S_1\circ\varTheta^{\varepsilon_1}\circ S_2 \circ\cdots\circ S_{\ell}\circ\varTheta^{\varepsilon_{\ell}})^{-1}$. Then, there are $z_1,w_1,z_2,\dots,w_{\ell} \in X_0$ so that  
\begin{equation*} 
x \overset{\varTheta^{-\varepsilon_0}}{\rightarrow} z_1 \overset{S_1^{-1}}{\rightarrow} w_1 \overset{\varTheta^{-\varepsilon_1}}{\rightarrow} z_2 \overset{S_2^{-1}}{\rightarrow} \cdots \overset{S_{\ell}^{-1}}{\rightarrow} w_{\ell} \overset{\varTheta^{-\varepsilon_{\ell}}}{\rightarrow} y. 
\end{equation*} 
Suppose here that $(x,y) \in \mathcal{S}$. Then, 
\begin{equation*} 
y \underset{\mathcal{S}}{\rightarrow} x \overset{\varTheta^{-\varepsilon_0}}{\rightarrow} z_1 \overset{S_1^{-1}}{\rightarrow} w_1 \overset{\varTheta^{-\varepsilon_1}}{\rightarrow} z_2 \overset{S_2^{-1}}{\rightarrow} \cdots \overset{S_{\ell}^{-1}}{\rightarrow} w_{\ell} \overset{\varTheta^{-\varepsilon_{\ell}}}{\rightarrow} y
\end{equation*} 
becomes a reduced word again, a contradiction. Therefore, the graph of $(\varTheta^{\varepsilon_0}\circ S_1\circ\varTheta^{\varepsilon_1}\circ S_2 \circ\cdots\circ S_{\ell}\circ\varTheta^{\varepsilon_{\ell}})^{-1}$ must be contained in $\mathcal{R}\setminus\mathcal{S}$ up to $\mu$-null set, and hence $E_{\mathcal{S}}^{\mathcal{R}}(w) = 0$.   
\end{proof} 

We have seen that any HNN equivalence relation is regarded as a particular case of HNN extensions of von Neumann algebras via Feldman-Moore's construction. Here, we would like to explain a close relation between the observation due to Gaboriau \cite[lines 12--26 in p.66]{Gaboriau:InventMath00} and ours (Proposition \ref{prop2.1}). Thanks to Proposition \ref{prop2.5} Gaboriau's observation is formulated in the framework of von Neumann algebras as follows. Let $N \supseteq D$ and $\theta : D \rightarrow N$ be as in Proposition \ref{prop2.1}, and suppose that $N$, $D$ and $\theta(D)$ all contain a common Cartan subalgebra, say $C$. By \cite{Aoi:JMathSocJapan03} there are unique (faithful normal) conditional expectations $E_D^N : N \rightarrow D$, $E_{\theta(D)}^N : N \rightarrow \theta(D)$, and we can consider the HNN extension $\big(M, E_N^M, u(\theta)\big) = \big(N,E_D^N\big)\bigstar_{D}\big(\theta,E_{\theta(D)}^N\big)$. Although there is no longer any reason supporting that $C$ becomes a MASA in $M$, any HNN equivalence relation gives such a triple thanks to Proposition \ref{prop2.5}. By Proposition \ref{prop2.1} $M\otimes M_2(\mathbf{C}) \supseteq C\otimes\mathbf{C}^2$ is isomorphic to $\mathcal{M} \supseteq \lambda(C\oplus C)$, by which one observes that $\lambda(C\oplus C)$ becomes a Cartan subalgebra in $\mathcal{M}$ when $C$ is an MASA in $M$. Here $\mathcal{M}$ is the amalgamated free product appeared in the construction of the HNN extension $M$ as in Proposition \ref{prop2.1}. This is nothing less than Gaboriau's observation in the von Neumann algebra context. Gaboriau's observation also consists of the converse assertion. Namely, he also stated, by giving an explicit description, that any amalgamated free product equivalence relation is stably isomorphic to a certain HNN equivalence relation. Its operator algebra counterpart will be explained briefly in the next subsection.     

\subsection{From Amalgamated free products to HNN extensions}  

Let $P_1, P_2, Q$ be $\sigma$-finite von Neumann algebras with two embeddings $\iota_1 : Q \hookrightarrow P_1$, $\iota_2 : Q \hookrightarrow P_2$. Suppose that there are two faithful normal conditional expectations $E_1 : P_1 \rightarrow \iota_1(Q)$, $E_2 : P_2 \rightarrow \iota_2(Q)$. Then, let $\big(P,E\big) := \big(P_1,E_1:\iota_1\big)\bigstar_Q\big(P_2,E_2:\iota_2\big)$ be the amalgamated free product of von Neumann algebras. Set $N := P_1\oplus P_2 \supseteq D := \iota_1(Q)\oplus\iota_2(Q)$, and define the bijective $*$-homomorphism $\theta : (\iota_1(x),\iota_2(y)) \in D \mapsto (\iota_1(y),\iota_2(x)) \in D$. Also define $E_D^N = E_{\theta(D)}^N := E_1\oplus E_2 : N \rightarrow D = \theta(D)$. Then, let $\big(M,E_N^M,u(\theta)\big) = \big(N,E_D^N\big)\bigstar_D\big(\theta,E_{\theta(D)}^N\big)$ be the HNN extension. Set $p := 1_{P_1}\oplus0 \in D$, and denote by $M_0$ the von Neumann subalgebra generated by $N$ and $v := pu(\theta)$ (a partial isometry with $v^* v = \theta(p) = 1-p$, $vv^* = p$). It is plain to see that $e_{11} := p$, $e_{12} := v$, $e_{21} := v^*$, $e_{22} := 1-p$ form a $2\times2$ matrix unit system in $M_0$, and moreover that $e_{11}M_0 e_{11}$ is generated by $e_{11}Ne_{11} = P_1\oplus0$ and $e_{12}N e_{21} = v(0\oplus P_2)v^* = u(\theta)(0\oplus P_2)u(\theta)$ (see e.g.~\cite[Lemma 5.2.1]{VoiculescuDykemaNica:Book}). The restriction $F := E_D^N\circ E_N^M\big|_{e_{11}M_0 e_{11}}$ clearly gives a faithful normal conditional expectation from $e_{11}M_0 e_{11}$ onto $e_{11}De_{11} = \iota_1(Q)\oplus0$. It is trivial that the restriction of $F$ to $e_{11}Ne_{11} = P_1\oplus0$ is given by $E_1\oplus0$. Also, the characterization of HNN extensions enables us to compute  
\begin{align*} 
F\big(u(\theta)(0\oplus x)u(\theta)^*\big) 
&= 
E_D^N\circ E_N^M\big(u(\theta)(0\oplus(E_2(x)+(x-E_2(x)))u(\theta)^*\big) \\
&= 
E_D^N\circ E_N^M\big(u(\theta)\theta(\iota_1(\iota_2^{-1}(E_2(x))\oplus0)u(\theta)^*\big) \\
&\phantom{aaaaaaaa}+ 
E_D^M\circ\underbrace{E_N^M\big(u(\theta)(0\oplus(x-E_2(x)))u(\theta)^*\big)}_{=0} \\
&= 
E_D^N\circ E_N^M\big(\iota_1\circ\iota_2^{-1}(E_2(x))\oplus0\big) \\
&= 
\iota_1\circ\iota_2^{-1}(E_2(x))    
\end{align*}
for $x \in P_2$. Define $\lambda : x \in Q \hookrightarrow \iota_1(x)\oplus0 \in e_{11}De_{11} \subseteq e_{11}Me_{11}$, $\lambda_1 : x \in P_1 \hookrightarrow x\oplus0 \in e_{11}Ne_{11} \subseteq e_{11}Me_{11}$, $\lambda_2 : x \in P_2 \hookrightarrow u(\theta)(0\oplus x)u(\theta)^* \in e_{12}Ne_{21} \subseteq e_{11}Me_{11}$. Then, we have  
\begin{gather*}
\lambda_1\circ\iota_1(x) =  \iota_1(x)\oplus0 = \lambda(x), \\ 
\lambda_2\circ\iota_2(x) = u(\theta)^*(0\oplus\iota_2(x))u(\theta)^* = u(\theta)\theta(\iota_1(x)\oplus0)u(\theta)^* = \iota_1(x)\oplus0 = \lambda(x)
\end{gather*} 
 for $x \in Q$. Since 
\begin{gather*} 
\mathrm{Ker}F\cap(P_1\oplus0) = \mathrm{Ker}E_1\oplus0 \subseteq \mathrm{Ker}E_D^N, \\ 
\mathrm{Ker}F\cap u(\theta)(0\oplus P_2)u(\theta)^* = u(\theta)(0\oplus\mathrm{Ker}E_2)u(\theta)^* 
\subseteq u(\theta)\mathrm{Ker}E_{\theta(D)}^N u(\theta)^*
\end{gather*} 
one easily derives, from the condition (M) in \S\S1.1, that $\lambda_1(P_1) = P_1\oplus0$ and $\lambda_2(P_2) = u(\theta)(0\oplus P_2)u(\theta)^*$ are free with respect to $F$.

Summarizing the discussion so far we conclude:  

\begin{proposition}\label{prop2.6} Let $\big(M,E_N^M,u(\theta)\big) = \big(N,E_D^N\big)\bigstar_D\big(\theta,E_{\theta(D)}^N\big)$ be the HNN extension with 
\begin{align*} 
&N := P_1\oplus P_2 \supseteq D := \iota_1(Q)\oplus\iota_2(Q), \\
&\theta : (\iota_1(x),\iota_2(y)) \in D \mapsto (\iota_1(y),\iota_2(x)) \in D, \\
&E_D^N = E_{\theta(D)}^N := E_1\oplus E_2, \\
&p\,(=e_{11}) := 1_{P_1}\oplus0 \in N,   
\end{align*}
and $M_0$ be the von Neumann subalgebra of $M$ generated by $N$ and $v := pu(\theta)$. Then,  
the compressed system $\big(pM_0 p,F = E_D^N\circ E_N^M\big|_{pM_0 p}\big)$ is identified with the amalgamated free product $(P,E) = (P_1,E_1:\iota_1)\bigstar_Q(P_2,E_2:\iota_2)$. 
\end{proposition} 

The subalgebra $M_0$, the conditional expectation $E_N^{M_0} := E_N^M|_{M_0}$ and the partial isometry $v$ can be charactereized, similarly as in the case of $(M,E_D^M,u(\theta))$, by the following two conditions: (A) $v\theta(d)v^* = d$ for every $d \in pD$; (M) $E_N^{M_0}(w)=0$ for every nonzero word $w = v^{\varepsilon_0}n_1 v^{\varepsilon_1}n_2\cdots n_{\ell}v^{\varepsilon_{\ell}}$ (with $n_1,\dots,n_{\ell} \in N$, $\varepsilon_0,\dots,\varepsilon_{\ell} \in \{\cdot,*\}$) which satisfies that $\varepsilon_{j-1}\neq\varepsilon_j$ implies that
\begin{itemize}
\item[] $n_j \in \mathrm{Ker}\big(E_{\theta(D)}^N|_{\theta(p)N\theta(p)}\big)$ 
when $\varepsilon_{j-1} = \cdot, \ \varepsilon_j = *$; 
\item[] $n_j \in \mathrm{Ker}\big(E_D^N|_{pNp}\big)$   
when $\varepsilon_{j-1}  = *, \ \varepsilon_j = \cdot$.  
\end{itemize}
Hence the triple $(M_0,E_N^{M_0},v)$ depends only on $\theta|_{pD}$ and $E_D^N|_{pNp}, E_{\theta(D)}^N|_{\theta(p)N\theta(p)}$ so that it should be called the ``(generalized) HNN extension by the partial $*$-isomorphism $\theta|_{pD} : pD \rightarrow \theta(p)N\theta(p)$ with respect to $E_D^N|_{pNp}, E_{\theta(D)}^N|_{\theta(p)N\theta(p)}$," whose details will be discussed elsewhere. Here a partial $*$-isomorphism means an injective unital $*$-homomorphism from a subalegebra whose unit is different from a given algebra into a compressed algebra of the given one. We should also remark that the same assertion as Proposition \ref{prop2.1} still holds true for $(M_0,E_N^{M_0},v)$. Namely, $M_0\otimes M_2(\mathbf{C})$ can be identified with the amalgamated free product 
\begin{equation*} 
\left(\begin{bmatrix} N & N \\ N & N \end{bmatrix}, \begin{bmatrix} E_D^N & 0 \\ 0 & E_{\theta(D)}^N \end{bmatrix} : \begin{bmatrix} \mathrm{Id}_D & \\ & \theta \end{bmatrix}\right) \underset{D\oplus D}{\bigstar} \left(\begin{bmatrix} D & pD \\ pD & D \end{bmatrix}, \begin{bmatrix} \mathrm{Id}_D & 0 \\ 0 & \mathrm{Id}_D \end{bmatrix} : \begin{bmatrix} \mathrm{Id}_D & \\ & \mathrm{Id}_D \end{bmatrix}\right)
\end{equation*}
in the same way as in Proposition \ref{prop2.1}.  

The same facts as above (including Proposition \ref{prop2.6}) is still valid in the $C^*$-algebra settings. The reduced setting is treated by the exactly same argument, but the universal setting needs to use the universality similarly to Proposition \ref{prop2.4}. In the course of the proof, one easily observes the following general fact: 

\begin{fact}\label{fact2.7} Let $B \supseteq C$ be unital $C^*$-algebras,  $\theta : C \rightarrow B$ be an injective unital $*$-homomorphism, and $p$ be a {\rm (}non-zero{\rm )} central projection in $C$. Write $C_0 := pC$ and $\theta_0 := \theta|_{C_0} : C _0 \rightarrow \theta(p)B\theta(p)$. Let  
\begin{equation*} 
\mathcal{A}_0 := \left(\begin{bmatrix} B & B \\ B & B \end{bmatrix} : \begin{bmatrix} \mathrm{Id}_C & \\ & \theta \end{bmatrix}\right) \underset{C\oplus C}{\bigstar^{\mathrm{univ}}} \left(\begin{bmatrix} C & C_0 \\ C_0 & C \end{bmatrix} : \begin{bmatrix} \mathrm{Id}_C & \\ & \mathrm{Id}_C \end{bmatrix}\right) 
\end{equation*}   
be the universal amalgamated free product of $C^*$-algebras with the canonical embedding maps $j$ {\rm (}into the amalgamated subalgebra{\rm )}, $j_{\theta}$ {\rm (}into the first free component{\rm )}, $j_1$ {\rm (}into the second free component{\rm )}. Then, the $C^*$-subalgebra $A_0$ {\rm (}inside the compressed algebra of $\mathcal{A}_0$ by $j\left(\begin{bmatrix} 1 & 0 \\ 0 & 0 \end{bmatrix}\right)${\rm )} generated by $j_{\theta}\left(\begin{bmatrix} b & 0 \\ 0 & 0 \end{bmatrix}\right)$, $b \in B$, and $v := j_1\left(\begin{bmatrix} 0 & p \\ 0 & 0 \end{bmatrix}\right) j_{\theta}\left(\begin{bmatrix} 0 & 0 \\ 1 & 0 \end{bmatrix}\right)$ is universal with subject to the algebraic equations $v\theta_0(c)v^* = c$ for all $c \in C_0$. Moreover, $A_0\otimes M_2(\mathbf{C})$ is identified with $\mathcal{A}_0$ by the same way as in Proposition \ref{prop2.4}.      
\end{fact}   

Hence the matrix trick we employ also provides the precise construction of ``universal HNN extensions by partial $*$-isomorphisms" (compare with the comment after Proposition \ref{prop2.4}).   

\section{Results} 

\subsection{Factoriality and Type classification}

\subsubsection{General Results} Let $N \supseteq D$ be $\sigma$-finite von Neumann algebras with an injective normal unital $*$-homomorphism $\theta : D \rightarrow N$, and then two faithful normal conditional expectations $E_D^N : N \rightarrow D$, $E_{\theta(D)}^N : N \rightarrow \theta(D)$ are given.  

\begin{assumption}\label{assumption3.1} Assume that there are two unitaries $v_1, v_{\theta} \in N$ and two faithful normal states $\varphi_1, \varphi_{\theta}$ on $D$ such that 
\begin{itemize} 
\item[(a)] $E_D^N\left(v_1^m\right) = E_{\theta(D)}^N\left(v_{\theta}^m\right) = 0$ as long as $m \neq 0$; 
\item[(b)] $v_1 \in N_{\varphi_1\circ E_D^N}$ and $v_{\theta} \in N_{\varphi_{\theta}\circ\theta^{-1}\circ E_{\theta(D)}^N}$. 
\end{itemize} 
\end{assumption} 

In what follows, we will use the notational rule in \S\S2.1. Namely, $\big(M, E_N^M, u(\theta)\big)$ is the HNN extension of $N$ by $\theta$ with respect to $E_D^N$ and $E_{\theta(D)}^N$, and also $\mathcal{M}$ is the associated amalgamated free product so that $\mathcal{M} \cong M\otimes M_2(\mathbf{C})$. In what follows, we use the usual notations for ultraproducts of von Neumann algebras. Namely, for a von Neumann algebra $L$ and a free ultrafilter $\omega \in \beta(\mathbb{N})\setminus\mathbb{N}$, $L^{\omega}$ denotes the ultraproduct of $L$ with respect to $\omega$. If a von Neumann subalgebra $K \subseteq L$ is the range of a faithful normal conditional expectation from $L$, then $K^{\omega}$ can be naturally regarded as a von Neumann subalgebra of $L^{\omega}$. Moreover, for a bijective normal $*$-homomorphism $\alpha : L_1 \rightarrow L_2$ between von Neumann algebras gives a unique bijective normal $*$-homomorphism $\alpha^{\omega} : L_1^{\omega} \rightarrow L_2^{\omega}$.    

\begin{proposition}\label{prop3.1} Under Assumption \ref{assumption3.1} we have 
\begin{equation*} 
\left\{\begin{bmatrix} v_1 & 0 \\ 0 & v_{\theta} \end{bmatrix}, \begin{bmatrix} 0 & u(\theta) \\ u(\theta)^* & 0 \end{bmatrix}\right\}'\cap \left(M\otimes M_2(\mathbf{C})\right)^{\omega} \subseteq \begin{bmatrix} D & 0 \\ 0 & \theta(D) \end{bmatrix}^{\omega}.  
\end{equation*} 
In particular, 
\begin{equation} 
\left(M\otimes M_2(\mathbf{C})\right)' \cap \left(M\otimes M_2(\mathbf{C})\right)^{\omega} = \left(M\otimes M_2(\mathbf{C})\right)' \cap \begin{bmatrix} D & 0 \\ 0 & \theta(D) \end{bmatrix}^{\omega}. \label{eq5}
\end{equation}  
\end{proposition}  
\begin{proof} Via the bijective $*$-homomorphism $\Phi$ in Proposition \ref{prop2.1} 
\begin{equation*} 
\left(M\otimes M_2(\mathbf{C}) \supseteq N\otimes M_2(\mathbf{C}), E^M_N\otimes\mathrm{Id}\right)
\end{equation*} 
is identified with 
\begin{equation*} 
\left(\mathcal{M}\supseteq\lambda_{\theta}\left(N\otimes M_2(\mathbf{C})\right), \mathcal{E}_{\theta}\right)
\end{equation*}
and correspondingly 
\begin{equation*} 
\begin{bmatrix} v_1 & 0 \\ 0 & v_{\theta} \end{bmatrix},\, \begin{bmatrix} 0 & u(\theta) \\ u(\theta)^* & 0 \end{bmatrix}\quad\text{with}\quad 
V:=\lambda_{\theta}\left(\begin{bmatrix} v_1 & 0 \\ 0 & v_{\theta} \end{bmatrix}\right),\,  
W:=\lambda_1\left(\begin{bmatrix} 0 & 1 \\ 1 & 0 \end{bmatrix}\right), 
\end{equation*} 
respectively. Hence, it suffices to show that 
\begin{equation} 
\left\{ V, W \right\}' \cap \mathcal{M}^{\omega} \subseteq \lambda\left(D\oplus D\right) ^{\omega} = \lambda_{\theta}\left(\iota_{\theta}\left(D\oplus D\right)\right)^{\omega}. \label{eq6}  
\end{equation} 
With letting $\psi\left((d_1, d_2)\right) := \frac{1}{2}\left(\varphi_1(d_1) + \varphi_{\theta}(d_2)\right)$, a faithful normal state on $D\oplus D$, Assumption \ref{assumption3.1} (b) implies that   
\begin{equation*} 
\sigma_t^{\psi\circ\iota_{\theta}^{-1}\circ E_{\theta}}\left(V\right) = V 
\end{equation*} 
for $t \in \mathbf{R}$, and hence $V \in \left(N\otimes M_2(\mathbf{C})\right)_{\psi\circ\iota_{\theta}^{-1}\circ E_{\theta}}$. Since 
\begin{equation*}
E_{\theta}\left(\begin{bmatrix} v_1 & 0 \\ 0 & v_{\theta}\end{bmatrix}^m\right) = 0\ \ (m\neq 0),\quad  
E_1\left(\begin{bmatrix} 0 & 1 \\ 1 & 0 \end{bmatrix}\right) = 0 
\end{equation*} 
thanks to Assumption \ref{assumption3.1} (a), we can apply \cite[Proposition 5]{Ueda:TAMS03} (note that the assumption ``$uDu^*=D=wDw^*$" there is never used in the proof as remarked in \cite[p.400]{Ueda:JFA05} so that we can apply it) to
\begin{equation*}
\left(\mathcal{M},\mathcal{E}\right) = \left(N\otimes M_2(\mathbf{C}), E_{\theta} : \iota_{\theta}\right) \underset{D\oplus D}{\bigstar}\left(D\otimes M_2(\mathbf{C}), E_1 : \iota_1\right)
\end{equation*} 
with $V, W$, and thus for all $X \in \left\{V\right\}' \cap \mathcal{M}^{\omega}$ we get 
\begin{equation*} 
\left\Vert W\left(X - \mathcal{E}^{\omega}(X)\right)\right\Vert_{\left(\psi\circ\lambda^{-1}\circ\mathcal{E}\right)^{\omega}} \leq 
\left\Vert WX - XW \right\Vert_{\left(\psi\circ\lambda^{-1}\circ\mathcal{E}\right)^{\omega}}. 
\end{equation*}
This inequality immediately implies \eqref{eq6}. 
\end{proof}  

Here is a simple (probably well-known) lemma needed for the derivation of a result on HNN extensions from Proposition \ref{prop3.1}. 

\begin{lemma}\label{lem3.2} {\rm (}e.g.~\cite[Lemma 2.1]{Popa:InventMath81}{\rm )} Let $P \supseteq Q$ be von Neumann algebras and $e \in Q$ be a projection. Then, $\left(eQe\right)'\cap ePe = Q'e\cap ePe = \left(Q'\cap P\right)e$. 
\end{lemma}  

\begin{theorem}\label{S3-T1} Under Assumption \ref{assumption3.1} we have 
\begin{align}
\mathcal{Z}(M) &= \left\{ x \in D\cap\theta(D)\cap N' : \theta(x) = x \right\}, \label{eq7} \\ 
M' \cap M^{\omega} &= \left\{ x \in D^{\omega}\cap\theta^{\omega}\left(D^{\omega}\right)\cap N' :  \theta(x) = x \right\}. \label{eq8}
\end{align}
Moreover, the core $\widetilde{M}$ satisfies that 
\begin{equation} 
\mathcal{Z}\big(\widetilde{M}\big) = \big\{ x \in \widetilde{D}\cap\widetilde{\theta}\left(\widetilde{D}\right)\cap\widetilde{N}' : \widetilde{\theta}(x) = x \big\}, \label{eq9} 
\end{equation} 
where we use the notations in \S\S1.2.  
\end{theorem} 
\begin{proof} Applying Lemma \ref{lem3.2} to \eqref{eq5} with $e = \begin{bmatrix} 1 & 0 \\ 0 & 0 \end{bmatrix}$ and $\begin{bmatrix} 0 & 0 \\ 0 & 1 \end{bmatrix}$ we get, respectively,    
\begin{equation*} 
M'\cap M^{\omega} = M'\cap D^{\omega}, \quad M'\cap M^{\omega} = M'\cap\theta(D)^{\omega} = M'\cap\theta^{\omega}\big(D^{\omega}\big).  
\end{equation*}  
The desired assertions immediately follow from these equations since $M$ is generated by $N$ and $u(\theta)$ and also $u(\theta)\theta^{\omega}(x)u(\theta)^* = x$ for all $x \in D^{\omega}$. (For more details we refer to  \cite[p.406--409]{Ueda:JFA05}.)    
\end{proof} 

In the next remark, we use the notations in \S\S1.2. 

\begin{remark}\label{rem3.4} The dual action $\big\{\vartheta^M_t\big\}_{t \in \mathbf{R}}$ on $\widetilde{M}$ is defined in such a way that $\vartheta^M_t\big|_M = \mathrm{Id}_M$ and $\vartheta^M_t(\lambda(s)) = e^{-its}\lambda(s)$ for $s,t \in \mathbf{R}$. Then, $\vartheta^M_t$ commutes with $\widetilde{\theta}$ for every $t \in \mathbf{R}$. In particular, \eqref{eq9} implies \eqref{eq7}.
\end{remark}
\begin{proof} The commutativity between $\vartheta^M_t$ and $\widetilde{\theta}$ is clear from their definitions. If \eqref{eq9} was true, then it would follow that  
\begin{align*} 
\mathcal{Z}(M) 
&= \mathcal{Z}\big(\widetilde{M}\big)^{\vartheta^M} \\
&= \big\{ x \in \widetilde{D}\cap\widetilde{\theta}\big(\widetilde{D}\big)\cap\widetilde{N}' : \widetilde{\theta}(x) = x,\, \vartheta^M_t(x) = x\, (t \in \mathbf{R})  \big\} \\
&= \big\{ x \in D\cap\theta(D)\cap N' : \theta(x) = x \big\}. 
\end{align*} 
Here we use \cite[Theorem XII. 1.1]{Takesaki:Book2} for $\widetilde{M}$ and $\widetilde{D}$ twice. Note that $D\cap\theta(D)\cap N' = D\cap\theta(D)\cap\widetilde{N}'$ thanks to the fact that $\mathrm{Ad}\lambda(s)$ acts on the center $\mathcal{Z}(N)$ trivially for every $s\in\mathbf{R}$.  
\end{proof} 

\subsubsection{The Cartan Subalgebra Case} Here, we consider and discuss a particular case; both $D$ and $\theta(D)$ are assumed to be Cartan subalgebras in $N$. Since any MASA in a von Neumann algebra contains its center, we note that the domains of $\theta$ and $\widetilde{\theta}$ must contain $\mathcal{Z}(N)$ and $\mathcal{Z}\big(\widetilde{N}\big)$, respectively. 

\begin{theorem}\label{thm3.5} If $N$ has no type I direct summand, then 
\begin{align} 
\mathcal{Z}(M) &= \{ x \in \mathcal{Z}(N) : \theta(x) = x \}, \label{eq10} \\ 
\mathcal{Z}\big(\widetilde{M}\big) &= \{ x \in \mathcal{Z}\big(\widetilde{N}\big) : \widetilde{\theta}(x) = x \}. \label{eq11}
\end{align}
Moreover, if $N$ is either of type II or a non-type I factor, then 
\begin{equation} 
M' \cap M^{\omega} = \left\{ x \in D^{\omega}\cap\theta^{\omega}\left(D^{\omega}\right)\cap N' :  \theta^{\omega}(x) = x \right\}. \label{eq12}
\end{equation}  
\end{theorem} 
\begin{proof} Since the core $\widetilde{N}$ is of type II (under the hypothesis of the first assertion), i.e., a direct sum of von Neumann algebras of type II$_1$ and type II$_{\infty}$, the argument of \cite[Lemma 4.2]{Ueda:Pacific99} enables us to confirm that Assumption \ref{assumption3.1} holds for $\widetilde{M} = \widetilde{N}\bigstar_{\widetilde{D}}\widetilde{\theta}$, and hence Theorem \ref{S3-T1} with the aid of Remark \ref{rem3.4} shows \eqref{eq10} and \eqref{eq11} since $D\cap\theta(D)\cap N' = \mathcal{Z}(N)$ and $\widetilde{D}\cap\widetilde{\theta}\big(\widetilde{D}\big)\cap\widetilde{N}' = \mathcal{Z}\big(\widetilde{N}\big)$. The last assertion is also shown similarly by combining Theorem \ref{S3-T1} with the argument of \cite[Lemma 4.2]{Ueda:Pacific99}.    
\end{proof} 

\begin{remarks}\label{rem3.6} Theorem \ref{thm3.5} implies the following facts{\rm :} 
\begin{itemize} 
\item[(1)] If $N$ is a non-type I factor, then so is $M$ thanks to \eqref{eq10}. 
\item[(2)] If $N$ is a type III$_1$ factor, then so is $M$ thanks to \eqref{eq11}. 
\item[(3)] When $N$ is a non-type I factor, if $M$ is a factor of type III$_0$ then so must be $N$ thanks to \eqref{eq11}. 
\item[(4)] If $N$ is a non-type I factor, then $M_{\omega} = M'\cap M^{\omega} \subseteq  D^{\omega}$ thanks to \eqref{eq12}. In this case, by the argument given in \cite[Theorem 8]{Ueda:TAMS03} $M$ is shown never to be strongly stable, i.e., $M \not\cong M\otimes R$ with the hyperfinite II$_1$ factor $R$.    
\end{itemize} 
\end{remarks}

\begin{theorem}\label{thm3.7} If $N$ is a factor of type II$_1$ or of type III$_{\lambda}$ with $\lambda \neq 0$, then there is a faithful normal state $\varphi$ on $D$ with $\big(N_{\varphi\circ E_D^N}\big)'\cap N = \mathbf{C}1$ {\rm (}$\varphi\circ E_D^N$ should be the unique tracial state in the type II$_1$ case{\rm )}, and moreover
\begin{equation} 
T(M) = \left\{ t \in T(N) : \big[D\varphi\circ\theta^{-1}\circ E_{\theta(D)}^N : D\varphi\circ E_D^N\big]_t = 1 \right\}. \label{eq13}
\end{equation}    
\end{theorem} 
\begin{proof} The first part of assertion holds true thanks to \cite[Lemma 4.2]{Ueda:Pacific99}; more precisely, one can construct two faithful normal states $\varphi$ and $\varphi_{\theta}$ on $D$ in such a way that 
\begin{itemize} 
\item there are unitaries $v_1 \in N_{\varphi\circ E_D^N}$, $v_{\theta} \in N_{\varphi_{\theta}\circ\theta^{-1}\circ E_{\theta(D)}^N}$ with $E_D^N\big(v_1^m\big) = 0$ and $E_{\theta(D)}^N\big(v_{\theta}^m\big) = 0$ as long as $m\neq0$; 
\item $\big(N_{\varphi\circ E_D^N}\big)'\cap N = \mathbf{C}1$ and $\big(N_{\varphi_{\theta}\circ\theta^{-1}\circ E_D^N}\big)'\cap N = \mathbf{C}1$. 
\end{itemize}
Let us define the faithful normal conditional expectation $\mathcal{E} : M\otimes M_2(\mathbf{C}) \rightarrow D \oplus \theta(D) = \begin{bmatrix} D & \\ & \theta(D) \end{bmatrix}$ and the faithful state $\psi$ on $D\oplus\theta(D)$ by 
\begin{align*} 
\mathcal{E}\left(\begin{bmatrix} m_{11} & m_{12} \\ m_{21} & m_{22} \end{bmatrix}\right) 
&:= \begin{bmatrix} E_D^N\circ E_N^M(m_{11}) & \\ & E_{\theta(D)}^M\circ E_N^M(m_{22}) \end{bmatrix}, \\ 
\psi\left(\begin{bmatrix} d_{11} & \\ & \theta(d_{22}) \end{bmatrix}\right) 
&:= \frac{1}{2}\left(\varphi(d_{11}) + \varphi_{\theta}(d_{22})\right). 
\end{align*}  
Clearly $V := \begin{bmatrix} v_1 & 0 \\ 0 & v_{\theta} \end{bmatrix}$ is in the centralizer $\left(M\otimes M_2(\mathbf{C})\right)_{\psi\circ\mathcal{E}}$, and the proof of Proposition \ref{prop3.1} shows that all $X \in \left\{ V \right\}' \cap \left(M\otimes M_2(\mathbf{C})\right)$ and $W_1, W_2 \in \mathrm{Ker}\mathcal{E}$ must satisfy that  
\begin{equation} 
\left\Vert W_1\left(X-\mathcal{E}(X)\right) \right\Vert_{\psi\circ\mathcal{E}} \leq \left\Vert W_1 X - X W_2\right\Vert_{\psi\circ\mathcal{E}}. \label{eq14} 
\end{equation}   
Let $t_0$ be a real number such that $\sigma_{t_0}^{\psi\circ\mathcal{E}} = \mathrm{Ad}U$ for some unitary $U \in M\otimes M_2(\mathbf{C})$, and set $W := \begin{bmatrix} 0 & 1 \\ 1 & 0 \end{bmatrix}$. Since $\sigma_t^{\psi\circ\mathcal{E}}(W) \in \mathrm{Ker}\mathcal{E}$, \eqref{eq14} shows that 
\begin{equation*} 
\left\Vert \sigma_{t_0}^{\psi\circ\mathcal{E}}(W)\left(U-\mathcal{E}(U)\right)\right\Vert_{\psi\circ\mathcal{E}} \leq 
\left\Vert \sigma_{t_0}^{\psi\circ\mathcal{E}}(W)U - UW\right\Vert_{\psi\circ\mathcal{E}} = 0. 
\end{equation*} 
Hence, $U=\mathcal{E}(U) = \begin{bmatrix} u & 0 \\ 0 & u_{\theta} \end{bmatrix}$ for some unitaries $u \in D$, $u_{\theta} \in \theta(D)$. 
It is plain to see that  
\begin{equation*}
\sigma_t^{\psi\circ\mathcal{E}}\left(\begin{bmatrix} m_{11} & m_{12} \\ m_{21} & m_{22} \end{bmatrix}\right) 
= \begin{bmatrix} 
 \sigma_t^{\varphi_1\circ E_D^N\circ E_N^M}(m_{11}) &  \sigma_t^{\varphi_1\circ E_D^N\circ E_N^M}(m_{12})u_t \\ 
u_t^* \sigma_t^{\varphi_1\circ E_D^N\circ E_N^M}(m_{21}) & u_t^* \sigma_t^{\varphi_1\circ E_D^N\circ E_N^M}(m_{22})u_t \end{bmatrix} 
\end{equation*} 
for $m_{ij} \in M$, $i,j=1,2$, and $t \in \mathbf{R}$ with letting $u_t := \big[D\varphi_1\circ E_D^N : D\varphi_{\theta}\circ\theta^{-1}\circ E_{\theta(D)}^N\big]_t$. Thus, we see that $\sigma_{t_0}^{\varphi\circ E_N^M\circ E_N^M} = \mathrm{Ad}u$. Since $N_{\varphi\circ E_D^N}$ sits in $M_{\varphi\circ E_D^N\circ E_N^M}$, we have $u \in \big(N_{\varphi\circ E_D^N}\big)'\cap D \subseteq \big(N_{\varphi\circ E_D^N}\big)'\cap N = \mathbf{C}1$ so that $\sigma_{t_0}^{\varphi\circ E_D^N\circ E_N^M} = \mathrm{Id}$. Consequently, $t \in T(M)$ if and only if $\sigma_t^{\varphi\circ E_D^N\circ E_N^M}=\mathrm{Id}$, which is equivalent to that $t \in T(N)$ and $\sigma_t^{\varphi\circ E_D^N\circ E_N^M}\big(u(\theta)\big) = u(\theta)$ since $M = \left\{ N, u(\theta) \right\}''$ and $\sigma_t^{\varphi\circ E_D^N\circ E_N^M}\big|_N = \sigma_t^{\varphi\circ E_D^N}$. Hence, the desired assertion immediately follows thanks to \eqref{eq1}.  
\end{proof}  

When $N$ is a type II$_1$ factor, the T-set $T(M)$ can be described more explicitly as follows. Let $\tau$ be the unique tracial state on $N$. Since $\big(\tau|_D\big)\circ E_D^N = \tau = \big(\tau|_{\theta(D)}\big)\circ E_{\theta(D)}^N$ must hold, we have $\big[D\big(\tau|_D\big)\circ E_D^N : D\big(\tau|_D\big)\circ\theta^{-1}\circ E_{\theta(D)}^N\big]_t = \big[D\big(\tau|_{\theta(D)}\big) : D\big(\tau|_D\big)\circ\theta^{-1}\big]_t$, and hence Theorem \ref{thm3.7} \eqref{eq13} is re-written as  
\begin{equation*} 
T(M) = \left\{ t \in \mathbf{R} :  \big[D\big(\tau|_{\theta(D)}\big) : D\big(\tau|_D\big)\circ\theta^{-1}\big]_t = 1\right\}. 
\end{equation*}

Assume that we have a type II$_1$ factor $N$ with two Cartan subalgebras $C_1$, $C_2$ and that $\tau_N$ is the unique tracial state on $N$.  Then, $C_1 \cong C_2 \cong L^{\infty}[0,1]$ in such a way that the Lebesgue measure $\nu$ on $[0,1]$ is the measure induced from the restrictions of $\tau_N$ to $C_1$ and $C_2$, respectively. Letting $D := C_1$ we can construct an bijective $*$-homomorphism $\theta : D \rightarrow N$ with $\theta(D) = C_2$. The $\theta$ induces the non-singular transformation $\varTheta$ on $[0,1]$. The above computation shows that 
\begin{equation*} 
T\big(N\bigstar_D\theta\big) = \Bigg\{ t \in \mathbf{R} :  \left(\Big(\frac{d\nu\circ\varTheta^{-1}}{d\nu}\Big)(\omega)\right)^{it} \equiv 1\ \text{for a.e.~$\omega\in[0,1]$} \Bigg\}.
\end{equation*}
Moreover one can see that the following are equivalent: (i) $N\bigstar_D\theta$ is semifinite. (ii) $N\bigstar_D\theta$ is a type II$_1$ factor. (iii) $\tau|_{\theta(D)} = (\tau|_D)\circ\theta^{-1}$.     

\subsection{A Sufficient Condition for Simplicity} Here, we will give a partial answer to the question of simplicity of reduced HNN extensions of $C^*$-algebras. Our method is to derive from a result on the simplicity of reduced amalgamated free products of $C^*$-algebras, due to K.~McClanahan \cite{McClanahan:CanJMath94} with the aid of Proposition 2.2. 

Let us first briefly review the above-mentioned result of McClanahan (which essentially comes from a technique due to D.~Avitzour \cite{Avitzour:TAMS82}). Let $P_1$, $P_2$, $Q$ be unital $C^*$-algebras and 
$\eta_1 : Q \hookrightarrow P_1$, $\eta_2 : Q \hookrightarrow P_2$ be embeddings. Assume that there are two conditional expectations $F_1 : P_1 \rightarrow \eta_2(Q)$, $F_2 : P_2 \rightarrow \eta_2(Q)$. Let $(P,F) := (P_1,F_1:\eta_1)\bigstar_{Q}(P_2,F_2:\eta_2)$ be the reduced amalgamated free product of $C^*$-algebras, where the canonical embeddings are denoted by $\rho : Q \hookrightarrow P$, $\rho_1 : P_1 \hookrightarrow P$, $\rho_2 : P_2 \hookrightarrow P$, which satisfies that $\rho = \rho_1\circ\eta_1 = \rho_2\circ\eta_2$ and $F : P \rightarrow \rho(Q)$ is a conditional expectation. Let us introduce the conditions:  
\begin{itemize} 
\item[1$^{\circ}$] There are unitaries $u,v \in P_1$, $w \in P_2$ such that 
%\begin{itemize} 
%\item 
$u\mathrm{Ker}F_1 u^* \subseteq \mathrm{Ker}F_1$, 
%\item 
$F_1(u^* v) = 0$, 
%\item 
$w\mathrm{Ker}F_2 w^* \subseteq \mathrm{Ker}F_2$; 
%\end{itemize} 
\item[2$^{\circ}$] For every $x \in Q$ and every $j \in \mathbb{Z}\setminus\{0\}$, there is an increasing sequence $\{m_k\}_{k=1,2,\dots}$ of natural numbers such that  
$$
\left[\rho(x), \left(\rho_1(u)\rho_2(w)\right)^{m_k}\rho_1(v)\rho_2(w)\rho_1(v)\left(\rho_2(w)\rho_1(u)\right)^j\right] = 0
$$
for all $k \geq k_0$ with some $k_0 \in \mathbb{N}$, 
\end{itemize}  
and then the subsets of $P_i$, $i=1,2$:  
\begin{equation*} 
\mathcal{N}^{(2)}(F_i) 
:= 
\left\{ (x,y) \in P_i\times P_i : x\mathrm{Ker}F_i y \subseteq \mathrm{Ker}F_i,\, x \eta_i(Q) y \subseteq \eta_i(Q) \right\}, 
\end{equation*} 
which act on $P$ by left-right multiplication. (Note that two more kinds of subsets are used in \cite{McClanahan:CanJMath94} to formulate the assertion, but they are nothing less than $Q$ and thus meaningless, since $Q$ is unital.) It is not so difficult to see that for any $(x,y) \in \mathcal{N}^{(2)}(F_i)$ one has $\rho_i(x) F(z)\rho_i(y) = F(\rho_i(x)z\rho_i(y))$ for every $z \in P$. Then, what McClanahan showed is: 

\begin{proposition}\label{prop3.8} {\rm (}\cite[Proposition 3.10]{McClanahan:CanJMath94}{\rm )} 
If the conditions  1$^{\circ}$, 2$^{\circ}$ hold, then any algebraic ideal $J \vartriangleleft P$ must satisfy that $F(J)$ sits inside the norm closure of $J$. Moreover, if $Q$ is further assumed to have no non-trivial $C^*$-ideal invariant under the actions of $\mathcal{N}^{(2)}(F_i)$, $i=1,2$, then $P$ must be simple.    
\end{proposition} 

The next lemma is shown by a simple calculation. 

\begin{lemma}\label{lem3.9} Assume that the unitaries $u, v, w$ in the condition 1$^{\circ}$ satisfy that $u,v \in \eta_1(Q)'\cap P_1$ and that $w^2 = 1$, i.e., $w$ is a self-adjoint unitary, and moreover $w\eta_2(Q)w = \eta_2(Q)$. Then, the condition 2$^{\circ}$ automatically holds true with $m_k := 2k-j-1$, $k \geq \frac{j+2}{2}$. 
\end{lemma} 

Lemma \ref{lem3.9} apparently gives the following variant of Proposition \ref{prop3.8}: 

\begin{proposition}\label{prop3.10} Assume that there are unitaries $u, v \in \eta_1(Q)'\cap P_1$ and $w=w^* \in P_2$ such that
\begin{itemize} 
\item $u \mathrm{Ker}F_1 u^* \subseteq \mathrm{Ker}F_1$, 
\item $F_1(u^* v) = 0$, 
\item $w \mathrm{Ker}F_2 w \subseteq \mathrm{Ker}F_2$, 
\item $w\eta_2(Q)w = \eta_2(Q)$, 
\item $\eta_2(Q)$ has no non-trivial $C^*$-ideal under the actions of $\mathcal{N}^{(2)}(F_i)$, $i=1,2$. 
\end{itemize} 
Then, $P$ must be simple. 
\end{proposition}    

We are now in position to apply McClanahan's result to the case of reduced HNN extensions. In what follows, we keep the setting and notations in \S\S2.2; namely, $\big(A, E_B^A,u(\theta)\big) = \big(B,E_C^B\big)\bigstar_C \big(\theta,E_{\theta(C)}^B\big)$ is a reduced HNN extension of $C^*$-algebras. Let  $\mathcal{N}^{(2)}\big(E_C^B\big), \mathcal{N}^{(2)}\big(E_{\theta(C)}^B\big)$ be defined as before, and they act on $B$ by left-right multiplication. We apply Proposition \ref{prop3.10} to the associated reduced amalgamated free product $(\mathcal{A},\mathcal{E})$ with letting $Q := C\oplus C$, $P_1 := B\otimes M_2(\mathbf{C})$, $P_2 := C\otimes M_2(\mathbf{C})$, $\eta_1 := \iota_{\theta}$, $\eta_2 := \iota_1$, $F_1 := E_{\theta}$, $F_2 := E_1$, and $P := \mathcal{A}$, $F := \mathcal{E}$, $\rho := \lambda$, $\rho_1 := \lambda_{\theta}$, $\rho_2 := \lambda_1$, and then get the following proposition:  

\begin{proposition}\label{prop3.11} Assume that there are unitaries $a \in C'\cap B$, $b \in \theta(C)'\cap B$ such that $E_C^B(a) = E_{\theta(C)}^B(b) = 0$ and that either $a\mathrm{Ker}E_C^B a^* \subseteq \mathrm{Ker}E_C^B${\rm ;} or $b\mathrm{Ker}E_{\theta(C)}^B b^* \subseteq \mathrm{Ker}E_{\theta(C)}^B$ holds. If $C$ has no $C^*$-ideal invariant under the actions of $\mathcal{N}^{(2)}\big(E_C^B), \mathcal{N}^{(2)}\big(E_{\theta(C)}^B\big)$ {\rm (}by left-right multiplication{\rm )}, then $A$ must be simple.      
\end{proposition} 
\begin{proof} Since $\mathcal{A} \cong A\otimes M_2(\mathbf{C})$ thanks to Proposition \ref{prop2.2}, it suffices to show that $\mathcal{A}$ is simple. We use Proposition \ref{prop3.10}, and thus need to specify the unitaries $u,v,w$ there in our setting. By symmetry we may and do assume that $E_C^B(a) = E_{\theta(C)}^B(b) = 0$ and $a\mathrm{Ker}E_C^B a^* \subseteq \mathrm{Ker}E_C^B$. Then, it is clear that the unitaries  
\begin{equation*} 
u := \begin{bmatrix} a & 0 \\ 0 & 1 \end{bmatrix}, \quad v := \begin{bmatrix} 1 & 0 \\ 0 & b \end{bmatrix}, \quad w := \begin{bmatrix} 0 & 1 \\ 1 & 0 \end{bmatrix} 
\end{equation*} 
satisfy the first four conditions in Proposition \ref{prop3.10}. Note that $(w,w) \in \mathcal{N}^2(F_2)$, and it is clear that any $C^*$-ideal in $Q = C\oplus C$ invariant under $\mathrm{Ad}w$ (via $\eta_1 = \iota_{\theta}$) must be of the form $C_0\oplus C_0$ with $C^*$-ideal $C_0 \vartriangleleft C$. 
%since $\mathrm{Ad}w$ switchs the first and the second diagonal components of $C\otimes\mathbf{C}^2 = \begin{bmatrix} C & \\ & C \end{bmatrix}$. 
Note also that $\mathcal{N}^{(2)}\big(E_C^B\big)\oplus\mathcal{N}^{(2)}\big(E_{\theta(C)}^B\big)$ are embedded into $\mathcal{N}^{(2)}(F_1)$ by 
\begin{gather*} 
\left((x_1,y_1),(x_2,y_2)\right) \in \mathcal{N}^{(2)}\big(E_C^B\big)\oplus\mathcal{N}^{(2)}\big(E_{\theta(C)}^B\big) 
%\phantom{aaaaaaaaaaaaaaaaaaaaaaaa}\\ \phantom{aaaaaaaaaaaaaaaaaaaaaaaa}
\mapsto \left(\begin{bmatrix} x_1  & \\ & x_2 \end{bmatrix},\begin{bmatrix} y_1  & \\ & y_2 \end{bmatrix}\right) \in \mathcal{N}^{(2)}(F_1), 
\end{gather*} 
respectively. Therefore, one easily observes that any $C^*$-ideal in $C\oplus C$ (considered inside $P$ via $\rho=\lambda$) that satisfies the hypothesis of Proposition \ref{prop3.10} must be of the form $C_0\oplus C_0$ with $C^*$-ideal $C_0 \vartriangleleft C$ invariant under the actions of $\mathcal{N}^{(2)}\big(E_C^B), \mathcal{N}^{(2)}\big(E_{\theta(C)}^B\big)$. By the assumption here, there is no such non-trivial $C^*$-ideal $C_0 \vartriangleleft C$, and hence $\mathcal{A}$ is simple by Proposition \ref{prop3.10}. 
\end{proof}  

\begin{example}\label{example3.12} {\rm Let $C$ be a simple $C^*$-algebra with a non-degenerate state $\varphi$. Set $B := C\otimes_{\mathrm{min}}C$, and identify the first component $C\otimes\mathbf{C}1$ in $B$ with $C$ itself. Then, we consider the injective unital $*$-homomorphism $\theta : x \in C = C\otimes_{\mathrm{min}}\mathbf{C}1 \mapsto 1\otimes x \in B$. The left and right slice maps of $\varphi$ give conditional expectations $E_C^B : B \rightarrow C$, $E_{\theta(C)}^B : B \rightarrow \theta(C)$, respectively. In this setting, if there is a unitary $u \in C$ such that $\varphi(u)=0$ and $\varphi\circ\mathrm{Ad}u = \varphi$, then the hypothesis of Proposition \ref{prop3.11} holds for the reduced HNN extension $\big(A,E_B^A,u(\theta)\big) = \big(B,E_C^B\big)\bigstar_C\big(\theta,E_{\theta(C)}^B\big)$, and hence $A$ is simple. }
\end{example}

Following \cite{delaHarpe:preprint05} we say a (discrete) group to be $C^*$-simple if its reduced group $C^*$-algebra is simple. The next corollary immediately follows from Proposition \ref{prop3.11}:  

\begin{corollary}\label{cor3.13} Let $G$ be a discrete group and $H$ be its subgroup with an injective homomorphism $\theta : H \rightarrow G$. If there are two elements $g_1 \in G\setminus H$ and $g_2 \in G\setminus\theta(H)$ so that $g_1$ and $g_2$ commute with $H$ and $\theta(H)$, respectively, and moreover $H$ is $C^*$-simple, then the HNN extension $G\bigstar_H\theta$ is $C^*$-simple.      
\end{corollary} 

The above corollary seems to be the first result on the $C^*$-simplicity of HNN extensions of groups. 

\subsection{$K$-Theory of HNN extensions} Our observation given in \S2 asserts that the computation of $K$-theory (also $KK$- and/or $E$-theory) of (universal and/or reduced) HNN extensions of $C^*$-algebras is reduced to that of the corresponding amalgamated free products. Here, we illustrate how to derive by obtaining the six terms exact sequence for $K$-groups associated with universal HNN extensions, which is exactly of the same kind of the one obtained in \cite{Anderson-Paschke:JOT86}.  

Here, we use (and keep) the setting and notations in \S\S2.3. Let us denote 
\begin{gather*} 
\mathcal{A}_1 := B\otimes M_2(\mathbf{C}), \quad 
\mathcal{A}_2 := C\otimes M_2(\mathbf{C}), \quad
\mathcal{B} := C\oplus C.   
\end{gather*}
and also the embedding map from a $C^*$-algebra $X = C$ or $B$ to another $Y = B$ or $A=B\bigstar_C^{\mathrm{univ}}\theta$ by $\iota_{X\hookrightarrow Y}$. Under a certain mild condition on $\mathcal{A}_1 \overset{\iota_{\theta}}{\hookleftarrow} \mathcal{B} \overset{\iota_1}{\hookrightarrow} \mathcal{A}_2$ it is known that the six terms exact sequence 
\begin{equation}\label{eq15}
\begin{matrix} 
K_0(\mathcal{B}) 
& \overset{(\iota_{\theta}{}_*,\iota_1{}_*)}{\longrightarrow} 
& K_0(\mathcal{A}_1)\oplus K_0(\mathcal{A}_2) 
& \overset{j_{\theta}{}_* - j_1{}_*}{\longrightarrow}  
& K_0(\mathcal{A}) \\
\uparrow
& & & & \downarrow \\
K_1(\mathcal{A}_1) 
& \underset{j_{\theta}{}_* - j_1{}_*} {\longleftarrow} 
& K_1(\mathcal{A}_1)\oplus K_1(\mathcal{A}_2) 
& \underset{(\iota_{\theta}{}_*,\iota_1{}_*)}{\longleftarrow} 
& K_1(\mathcal{B})
\end{matrix}
\end{equation} 
holds true. In fact, the most general result of this type was provided by K.~Thomsen \cite{Thomsen:JFA03}, where he assumed that $\mathcal{B}$ is nuclear or the existence of conditional expectations from $\mathcal{A}_1$, $\mathcal{A}_2$ onto $\iota_{\theta}(\mathcal{B})$, $\iota_1(\mathcal{B})$, respectively. Note that his conditions are apparently translated in our setup to the nuclearity of $C$ or the existence of conditional expectations from $B$ onto $C$, $\theta(C)$. 

Notice here that we have the following isomorphisms: 
\begin{gather*} 
K_0(\mathcal{B}) \cong K_0(C)\oplus K_0(C) \quad \text{by} \quad [(p,q)] \leftrightarrow [p]\oplus[q]\,; \\
\begin{cases} K_0(\mathcal{A}_1) \cong K_0(B) \\ K_0(\mathcal{A}_2) \cong K_0(C) \end{cases} \quad 
\text{with} \quad \left[\begin{pmatrix} p & \\ & q \end{pmatrix}\right] \leftrightarrow [p]+[q]\,; \\
K_0(\mathcal{A}) {\cong} K_0(A\otimes M_2(\mathbf{C})) \quad \text{by} \quad \Phi_*\,; \\
K_0(A\otimes M_2(\mathbf{C})) \cong K_0(A) \quad 
\text{with} \quad \left[\begin{pmatrix} p & \\ & q \end{pmatrix}\right] \leftrightarrow [p]+[q].        
\end{gather*} 
For the description of the second and the fourth isomorphisms above we use the obvious identification $M_n(D\otimes M_2(\mathbb{C})) = M_2(M_n(D))$ with an arbitrary $C^*$-algebra $D$, which identifies $M_n(D\otimes\mathbf{C}^2)$ with the diagonal matrices whose entries are from $M_n(D)$. By these facts we can re-write the upper horizontal line in \eqref{eq15} as follows.   
\begin{equation*} 
K_0(C)\oplus K_0(C)  \overset{\phi_0}{\longrightarrow}  K_0(B)\oplus K_0(C) \overset{\psi_0}{\longrightarrow}  K_0(A), 
\end{equation*} 
where the left arrow is given by $\phi_0 : [p]\oplus[q] \mapsto \big([p]+\theta_*([q])\big)\oplus\big([p]+[q]\big)$ and the right one by $\psi_0  : [p]\oplus[q] \mapsto [p]-[q]$. Similarly we have 
\begin{gather*} 
K_1(\mathcal{B}) \cong K_1(C)\oplus K_1(C) \quad \text{by} \quad [(u,v)] \leftrightarrow [u]\oplus[v]\,; \\
\begin{cases} K_1(\mathcal{A}_1) \cong K_1(B) \\ K_1(\mathcal{A}_2) \cong K_1(C) \end{cases} \quad 
\text{with} \quad \left[\begin{pmatrix} u & \\ & v \end{pmatrix}\right] \leftrightarrow [u]+[v]\,; \\
K_1(\mathcal{A}) {\cong} K_1(A\otimes M_2(\mathbf{C})) \quad \text{by} \quad \Phi_*\,; \\
K_1(A\otimes M_2(\mathbf{C})) \cong K_1(A) \quad 
\text{with} \quad \left[\begin{pmatrix} u & \\ & v \end{pmatrix}\right] \leftrightarrow [u]+[v],     
\end{gather*} 
and thus the lower horizontal arrow in \eqref{eq15} is also re-written to be 
\begin{equation*} 
K_1(A) \underset{\psi_1}{\longleftarrow} K_1(B)\oplus K_1(C) \underset{\phi_1}{\longleftarrow} K_1(C)\oplus K_1(C),
\end{equation*}
where $\phi_1 : [u]\oplus [v] \mapsto \big([u]+\theta_*([v])\big)\oplus\big([u]+[v]\big)$ and $\psi_1 : [u]\oplus[v] \mapsto [u]-[v]$. Hence, the six terms exact sequence \eqref{eq15} becomes 
\begin{equation}\label{eq16}
\begin{matrix} 
K_0(C)\oplus K_0(C) & \overset{\phi_0}{\longrightarrow} & K_0(B)\oplus K_0(C) & \overset{\psi_0}{\longrightarrow} & K_0(A) \\
\uparrow & & & & \downarrow \\
K_1(A) & \underset{\psi_1}{\longleftarrow} & K_1(B)\oplus K_1(C) & \underset{\phi_1}{\longleftarrow} & K_1(C). 
\end{matrix} 
\end{equation}   
Let $\phi'_0$ be the projection map from $K_0(C)\oplus K_0(C)$ to the second component and set $\psi'_0 : = (\mathrm{id}_B)_* - (\iota_{C\hookrightarrow B})_*$. Then, we have  
\begin{equation*} 
\begin{matrix} 
K_0(C) 
& \overset{\theta_* - (\iota_{C\hookrightarrow B})_*}{\longrightarrow} 
& K_0(B) \\
\\
\phi'_0\, \uparrow\quad 
& \circlearrowleft 
& \quad\uparrow\, \psi'_0 & \ \circlearrowleft \quad\searrow 
& (\iota_{B\hookrightarrow A})_* \\
\\
K_0(C)\oplus K_0(C) 
& \underset{\phi_0}{\longrightarrow} 
& K_1(B)\oplus K_0(C) 
& \underset{\psi_0}{\longrightarrow} 
&K_0(A), \\
\end{matrix}
\end{equation*} 
and 
%\begin{equation*} 
$\phi'_0\left(\mathrm{Ker}\phi_0\right) = \mathrm{Ker}\left(\theta_*-(\iota_{C\hookrightarrow B})_*\right)$. 
%\end{equation*} 
Similarly, let $\phi'_1$ be the projection from $K_1(C)\oplus K_1(C)$ onto the second component and set $\psi'_1 := (\mathrm{id}_B)_* - (\iota_{C\hookrightarrow B})_*$. Then, we have 
\begin{equation*} 
\begin{matrix} 
K_1(A)
& \underset{\psi_1}{\longleftarrow} 
& K_1(B)\oplus K_1(C) 
& \underset{\phi_1}{\longleftarrow} 
&K_1(C)\oplus K_1(C) \\
\\
(\iota_{B\hookrightarrow A})_* 
& \nwarrow\quad \circlearrowleft\ 
& \psi'_1\,\downarrow\quad 
& \circlearrowleft 
& \quad\downarrow\,\phi'_1 \\
\\
& 
& K_1(B) 
& \underset{\theta_*-(\iota_{C\hookrightarrow B})_*}{\longleftarrow} 
& K_1(C),
\end{matrix}
\end{equation*} 
and 
%\begin{equation*} 
$\phi'_1\left(\mathrm{Ker}\phi_1\right) = \mathrm{Ker}\left(\theta_*-(\iota_{C\hookrightarrow B})_*\right)$. 
%\end{equation*} 
From these facts together with \eqref{eq16} we finally get the following: 

\begin{proposition}\label{prop3.14} If $C$ is nuclear or there are conditional expectations from $B$ onto $C$ and $\theta(C)$, then the universal HNN extension $A = B\bigstar^{\rm{univ}}_C\theta$ satisfies the following six terms exact sequence:  
\begin{equation}\label{eq17} 
\begin{matrix} 
K_0(C) & \overset{\theta_* - (\iota_{C\hookrightarrow B})_*}{\longrightarrow} & K_0(B) & \overset{(\iota_{B\hookrightarrow A})_*}{\longrightarrow} & K_0(A) \\ 
\uparrow & & & & \downarrow \\
K_1(A) & \underset{(\iota_{B\hookrightarrow A})_*}{\longleftarrow} & K_1(B) & \underset{\theta_*-(\iota_{C\hookrightarrow B})_*}{\longleftarrow} & K_1(C).
\end{matrix}
\end{equation}
\end{proposition} 

\begin{remark}\label{rem3.15}{\rm Note that the above proposition apparently includes the celebrated six terms exact sequence for crossed-products by the integers $\mathbb{Z}$ due to M.~Pimsner and D.~Voiculescu \cite{PimsnerVoiculescu:JOT82} as the special case where $B = C$ and $\theta \in \mathrm{Aut}(B)$, i.e, $B\bigstar^{\rm{univ}}_C\theta \cong B\rtimes_{\theta}\mathbb{Z}$. The work \cite{PimsnerVoiculescu:JOT82} also deals with crossed-products by free groups $\mathbb{F}_n$ whose universal construction version can be also treated in the same way. However, we need to generalize what we have done in this paper to the setup of HNN extensions $B\bigstar^{\rm{univ}}_C\Theta$ with families $\Theta$ of injective unital $*$-homomorphisms from $C$ into $B$, see the comment after Proposition \ref{prop2.4}.}
\end{remark}   

We emphasize that our method of getting \eqref{eq17} still works even for reduced HNN extensions when the initial six terms exact sequence \eqref{eq15} holds true for the associated reduced amalgamated free products. In this direction, the main future problem is apparently to establish the $K$-amenability, i.e., the natural surjective homomorphism between the $K$-groups of reduced and universal HNN extensions is injective, under suitable assumptions. The same question for amalgamated free products was discussed by E.~Germain \cite{Germain:FreeProbabilityTheory97} but is not yet settled at the present moment ({\it n.b.} it was already settled only in the case where the amalgamated subalgebra consists only of the scalars $\mathbf{C}$, see \cite{Germain:Duke96}).  

\appendix 

\section{More on Factoriality} 

Here, we prove a certain relative commutant property for HNN extensions of von Neumann algebras. It is proved in the same line as before based on \cite[Appendix I]{Ueda:Pacific99}, which comes from works due to Avitzour \cite{Avitzour:TAMS82} and McClanahan \cite{McClanahan:CanJMath94}. 

Keep the notational rule in \S\S2.1,  i.e., $\big(M,E_N^M,u(\theta)\big) = \big(N,E_D^N\big)\bigstar_D \big(\theta,E_{\theta(D)}^N\big)$ is an HNN extension of von Neumann algebras and $\big(\mathcal{M},\mathcal{E}\big) := \big(N\otimes M_2(\mathbf{C}),E_{\theta} : \iota_{\theta}\big)\bigstar_{D\oplus D}\big(D\otimes  M_2(\mathbf{C}),E_1 : \iota_1\big)$ is the associated amalgamated free product of von Neumann algebras with canonical embedding maps $\lambda, \lambda_{\theta}, \lambda_1$ of $D\oplus D, N\otimes M_2(\mathbf{C}), D\otimes M_2(\mathbf{C})$, respectively, into $\mathcal{M}$. In what follows, we assume the following conditions: 
\begin{itemize}
\item[1$^{\circ}$] There is a faithful normal state $\varphi$ on $D$ so that $\varphi\circ E_D^N = \varphi\circ\theta^{-1}\circ E_{\theta(D)}^N$, and we denote it by $\psi$. 
\item[2$^{\circ}$] There are unitaries $a, b \in N_{\psi}$ such that $E_D^N(a)=0$, $E_{\theta(D)}^N(b)=0$ and $E_D^N\circ\mathrm{Ad}a = \mathrm{Ad}a\circ E_D^N$, $E_{\theta(D)}^N\circ\mathrm{Ad}b = \mathrm{Ad}b\circ E_{\theta(D)}^N$. 
\end{itemize}  
We then write
\begin{equation*} 
u = \begin{bmatrix} a & 0 \\ 0 & 1 \end{bmatrix}, \quad 
v = \begin{bmatrix} 1 & 0 \\ 0 & b \end{bmatrix} \in N\otimes M_2(\mathbf{C}) \quad\text{and}\quad  
w = \begin{bmatrix} 0 & 1 \\ 1 & 0 \end{bmatrix} \in D\otimes M_2(\mathbf{C}),  
\end{equation*}  
and it is easy to see that these $u,v$ and $w$ are in the centralizers $\big(N\otimes M_2(\mathbf{C})\big)_{\psi\otimes\mathrm{tr}_2}$ and $\big(N\otimes M_2(\mathbf{C})\big)_{\varphi\otimes\mathrm{tr}_2}$, respectively, where $\mathrm{tr}_2$ denotes the normalized trace on $M_2(\mathbf{C})$. Also, by the condition 2$^{\circ}$ we have   
\begin{gather*} 
E_{\theta}\circ\mathrm{Ad}u = \mathrm{Ad}u\circ E_{\theta}, \quad 
E_{\theta}\circ\mathrm{Ad}v = \mathrm{Ad}v\circ E_{\theta}, \quad 
E_1\circ\mathrm{Ad}w = \mathrm{Ad}w\circ E_1, \\
E_{\theta}(u) = E_{\theta}(v) = E_1(w) = 0.  
\end{gather*} 
Define 
\begin{equation*} 
\tilde{\varphi}\left((d_1, d_2)\right) := \frac{1}{2}\varphi(d_1+d_2), \quad 
(d_1, d_2) \in D\oplus D, 
\end{equation*}  
and then by the assumption 1$^{\circ}$ we see that $\widetilde{\varphi}\circ E_{\theta} = \psi\otimes\mathrm{tr}_2$ and $\widetilde{\varphi}\circ E_1 = \varphi\otimes\mathrm{tr}_2$. Hence, we can use \cite[Proposition I-C]{Ueda:Pacific99} so that  
\begin{equation*} 
\left\{ \lambda_{\theta}(u), \lambda_{\theta}(v), \lambda_1(w) \right\}' \cap \mathcal{M} \subseteq \lambda(D\otimes M_2(\mathbf{C})). 
\end{equation*} 
Note here that $u(\theta) \in M_{\psi\circ E_N^M}$ thanks to the condition 1$^{\circ}$ together with \eqref{eq1}, and thus the next proposition follows in the exactly same way as Theorem \ref{S3-T1}. 

\begin{proposition}\label{propB.1} Under the conditions 1$^{\circ}$, 2$^{\circ}$ one has
\begin{equation*} 
\big(M_{\psi\circ E_N^M}\big)'\cap M \subseteq \left\{ x \in D \cap \theta(D) \cap \big(N_{\psi}\big)' : \theta(x) = x \right\}. 
\end{equation*} 
If it is further assumed that $\psi$ is a tracial state, then 
\begin{equation*} 
\mathcal{Z}(M) = \left\{ x \in D \cap \theta(D) \cap N' : \theta(x) = x \right\}. 
\end{equation*} 
\end{proposition} 

We do not know whether the first inclusion relation is actually the equality or not, because it is not obvious whether $M_{\psi\circ E_N^M}$ is generated by $N_{\psi}$ and $u(\theta)$.

\end{document}